\documentclass[11pt,a4paper]{article}
\UseRawInputEncoding

\usepackage[francais]{babel}
\usepackage{enumitem}
\usepackage[utf8]{inputenc}
\usepackage[T1]{fontenc}
\usepackage{amsmath}
\usepackage{amsthm}
\usepackage{amssymb}
\usepackage{graphicx}
\usepackage{caption}
\usepackage{mathrsfs}
\usepackage{dsfont}
\usepackage[backref]{hyperref}
\usepackage{xcolor}
\usepackage{tikz-cd}
\usepackage{cancel}
\usepackage{ulem}
\usepackage{comment}
\usepackage{enumitem}   

\hypersetup{colorlinks=true, linkcolor=blue!60!black, urlcolor=red!60!black}
\usepackage{mathtools}

\newtheorem{definition}{Definition}
\newtheorem{prop}{Proposition}
\newtheorem{lemme}{Lemma}
\newtheorem{theoreme}{Theorem}
\newtheorem*{theoremee}{Theorem}
\newtheorem{theorem}{Theorem}

\newtheorem{remarque}{Remark}
\newtheorem*{exemple}{Example}
\newtheorem*{exemples}{Examples}
\newtheorem{corollaire}{Corollary}

\newenvironment{preuve}{\begin{proof} \rm}{\end{proof}}

\newcommand{\C}{\mathbb{C}}
\newcommand{\R}{\mathbb{R}}

\newcommand{\F}{\mathcal{F}}

\newcommand{\TF}{T\mathcal{F}}
\newcommand{\nF}{\nu\mathcal{F}}
\newcommand{\AqNps}{\mathcal{A}^{0,l}(\mathcal{N}_k^*)}
\newcommand{\AqNss}{\mathcal{A}^{0,l}(\mathcal{N}_*^*)}
\newcommand{\AAqNps}{A^{0,l}(\mathcal{N}_k^*)}

\newcommand{\Apq}{\mathcal{A}^{k,l}_X}
\newcommand{\AApq}{A^{k,l}_X}
\newcommand{\AqTFps}{\mathcal{A}^{0,l}(\bigwedge^k T\mathcal{F}^*)}
\newcommand{\AqTFpsplus}{\mathcal{A}^{0,l}(\bigwedge^k T\mathcal{F}^{*+1})}
\newcommand{\AAqTFps}{A^{0,l}(\bigwedge^k T\mathcal{F}^*)}
\newcommand{\pF}{\partial_\mathcal{F}}
\newcommand{\AqLpTX}{\mathcal{A}^{0,l}_X\otimes\bigwedge^k TX}
\newcommand{\AAqLpTX}{A^{0,l}_X(\bigwedge^k TX)}
\newcommand{\AAqLpTF}{A^{0,l}(\bigwedge^{k}T\mathcal{F})}
\newcommand{\AqLpTFs}{\mathcal{A}^{0,l}(\bigwedge^{k}T\mathcal{F}^*)}
\newcommand{\AAqLpTFs}{A^{0,l}(\bigwedge^{k}T\mathcal{F}^*)}
\newcommand{\AAqLkmpTFs}{A^{0,l}(\bigwedge^{p-k}T\mathcal{F}^*)}

\title{Calabi-Yau Foliations and Deformations}
\author{R\'emi \bsc{Danain-Bertoncini} \\}
\date{2024}

\begin{document}

\renewcommand{\proofname}{Proof}

\maketitle

\renewcommand{\contentsname}{Summary} 

\tableofcontents

\newpage

\section*{Abstract}

We propose in this article the study of the deformations of a Calabi-Yau type foliations $\F$. For three different types of deformations (unfoldings, holomorphic, transversally holomorphic) there exist Kuranishi spaces $K^f,K^h,K^{tr}$ parametrizing the corresponding families of deformations. We show that $K^f$ is smooth, and that we can obtain $K^h$ as the product $K^f\times K^{tr}$. At last, we show that we can see the $f$-deformations of $\F$ as the $tr$-deformations of a supplementary foliation $\mathcal{G}$.

\section*{Introduction}

It is customary in deformation theory to seek conditions that ensure the smoothness of the Kuranishi space for a given object. This property is obtained by assuming, for example, that there are no obstructions to the deformation of the object under consideration, which is naturally done by assuming the cancellation of a cohomology group in which obstructions are constructed. It is well known that Calabi-Yau manifolds, even without this cohomological assumption, have a smooth Kuranishi space. 

Considering holomorphic foliations on complex manifolds, we are interested in knowing whether there is a Calabi-Yau type hypothesis that ensures the non-obstruction of deformations of these foliations. For a well-chosen definition of what a strongly Calabi-Yau is, we show :

\begin{theoremee}
    Let $(X,\mathcal{F})$ be a compact Kähler manifold foliated by a strongly Calabi-Yau foliation. Then, the Kuranishi space $K^f$ of its $f$-deformation is smooth. 
\end{theoremee}

Those $f$-deformations are one of the three type of deformations that appears naturally for holomorphic foliations. The two remaining are $h$-deformations and $tr$-deformations. In idea, the $f$-deformations govern the deformations of the complex structures of the leaves of $\F$ (but not exactly), the $tr$-deformations govern the deformations of the transverse structure of $\F$ (forgetting a priori the holomorphic structure of the leaves) and the $h$-deformations a mix of both. For those three type of deformations, the existence of a Kuranishi space in each case has been proved in \cite{GHS,GiNi}. Here, we relate together thoses spaces as :

\begin{theoremee}
    Let $K^f_\F$, $K^{tr}_\F$ and $K^h_\F$, the Kuranishi spaces corresponding respectively to the $f,tr,h$-deformations of a regular holomorphic foliation $(X,\mathcal{F})$ admitting a regular holomorphic foliation supplementary to $\mathcal{F}$ and such that $K^f_\F$ is smooth.

    Then there exists an isomorphism between analytic spaces $\alpha:K^f_\F\times K^{tr}_\F\to K^h_\F$ making commutative the following diagramm :
    \begin{equation*}
   \begin{tikzcd}
  K^f_\F\times K^{tr}_\F \arrow[rr, "\alpha"] \arrow[rd, "p_2"]
  & & K^h_\F \arrow[ld,"\pi" ]\\
  & K^{tr}_\F &  \\ 
\end{tikzcd} 
\end{equation*}
\end{theoremee}

Finally, since our definition of strongly Calabi-Yau implies automatically the existence of a supplementary (i.e everywhere transverse and complementary) foliation $\mathcal{G}$ for $\F$, we can try to find relations between their Kuranishi spaces. We can relate the $f$-deformations of one to the $tr$-deformations of the other. It is the content of the following :

\begin{theoremee}\label{KtrKf}
   Let $(X,\F,\mathcal{G})$ be two supplementary regular holomorphic foliations on $X$ a compact complex manifold. We have :
   \begin{equation*}
       K^{tr}_\F\simeq  K^f_\mathcal{G}.
   \end{equation*}
   In particular, for $(X,\F)$ a strongly Calabi-Yau foliation and $\mathcal{G}$ a foliation supplementary to it, we have $K^{tr}_\F\simeq K^f_\mathcal{G}$.
   If moreover $X$ is Calabi-Yau, we get that $K^f_\F$, $K^{tr}_\F$ and $K^h_\F$ are smooth.
\end{theoremee}

We firstly give some preliminaries about the structures of interest and different points of view that can be taken to deal with our problems. In the section \ref{unobstruction} we give the results we need to obtain the smoothness of $K^f$. Then, in the section \ref{dec}, we show the decomposition theorem thanks to the almost complex multifoliate description of the holomorphic foliations. At the end, we relate $K^{tr}_\F$ and $K^f_\mathcal{G}$.

\subsection*{Acknowledgements}

I warmly thank Christophe Mourougane for our regular discussions, Fr\'ed\'eric Touzet for the discussion of the two notions of Calabi-Yau foliations in this text, Jorge Vitório Pereira for the discussions and the welcome he gave me at IMPA in Rio de Janeiro, and Marcel Nicolau who kindly enlightened me on various results discussed in this text. (Avec le soutien financier du programme d'aides à la mobilit\'e
internationale op\'er\'e par le Collège doctoral de Bretagne et cofinanc\'e avec la R\'egion
Bretagne [et Rennes M\'etropole])

\newpage

\section{Preliminaries}

All over the text, $M$ will be a compact differential manifold and $X$ a compact complex manifold whose underlying differential manifold is $M$. Depending on the context, we will denote $T\F$ for the tangent bundle of $\F$ as a differential foliation but also for the tangent bundle of $\F$ as a holomorphic foliation. Sometime, to make it clear, we use the notation $T\F^{0,1}:=(T\F\otimes \C)\cap T^{1,0}$. Recall also that $\nu\F:=TM/T\F$ or $\nu\F:=TX/T\F$, when $\F$ seen as a differential respectively holomorphic foliation, depending on the context.

\subsection{$\Gamma$-structures and examples}\label{Gammastructures}

A natural way of defining certain geometric structures (of differential manifolds, complex manifolds, foliated manifolds, ...) on a topological space $S$ is to give local coordinates whose transition maps belong to a given pseudogroup. More precisely :

\begin{definition}
    A pseudogroup $G$ on a topological space $S$ is given by a collection of homeomorphisms $(g:U_g\to V_g)_{g\in G}$ between open sets of $S$ such that :
    \begin{enumerate}
        \item $(U_g)_g$ is a cover of $S$,  
        \item the restriction $g_{|U'}$ of an element $g\in G$ to any $U'\subset U_g$ is an element of $G$,
        \item the composition of two elements of $G$, where it is defined, is an element of $G$,
        \item the inverse of an element of $G$ is an element of $G$,
        \item if $h:U\to V$ is a homeomorphism between open sets of $S$ such that the restriction to any non-empty open set of the form $U\cap U_g$ is an element of $G$, then $h$ is an element of $G$.
    \end{enumerate}
\end{definition}

Therefore we define a $\Gamma$-structure as the datum of an atlas whose transition maps belong to a well-adapted pseudogroup :

\begin{definition}
    Let $\Gamma$ be a pseudogroup on a model topological space $\mathcal{M}$. A $\Gamma$-structure on a topological space $\mathcal{X}$ is the datum of an atlas $(U_i,\phi_i)_i$ consisting of open sets $U_i$ covering $\mathcal{X}$ and homeomorphisms between open sets of $\mathcal{X}$ and $\mathcal{M}$ such that the transition maps $\phi_i\circ(\phi_j)^{-1}$, where they are defined, are elements of $\Gamma$.
\end{definition}

The structures of interest in this text are respectively those associated with the foliation (its holomorphic structure or its transversely holomorphic structure) of codimension $q$ (and $p:=n-q$), and the structures supported by the families of deformations we consider, parametrized by an analytic space $R$. These structures are obtained as $\Gamma$-structures for which we specify the pseudogroups :\\

\begin{tabular}{|p{5cm}|p{8cm}|}
    \hline
    \textbf{Structures} & \textbf{Associated pseudogroup} \tabularnewline
    \hline
    Transversaly holomorphic foliation & $\Gamma^{tr}$ the pseudogroup of local automorphisms of $ \mathbb{R}^{2n-2q}\times \mathbb{C}^q$ of the form $F(x,y,w)=(f_{l}(x,y,x',y'),f_{tr}(w))$, where $f_{tr}$ depending holomorphically on $w:=x'+iy'$ and $f_{l}$ depending differntially on $x,y,x'$ and $y'$.  \tabularnewline
    \hline
    Holomorphic foliation & $\Gamma^h$ the pseudogroup of local automorphisms of $ \mathbb{C}^{p}\times \mathbb{C}^q$ of the form $F(z,w)=(f_{l}(z,w),f_{tr}(w))$, where $f_{tr}$ depending holomorphically on $w$ and $f_{l}$ depending holomorphically on $z$ and $w$. \tabularnewline
    \hline
    Holomorphic foliation admitting a supplementary holomorphic foliation  & $\Gamma^{supp}$ the pseudogroup of local automorphisms of $ \mathbb{C}^p\times \mathbb{C}^q$ of the form $F(z,w)=(f_{l}(z),f_{tr}(w))$, where $f_{tr}$ depending holomorphically on $w$ and $f_{l}$ depending holomorphically on $z$.\tabularnewline
    \hline
    Family of $tr$-deformations & $\Gamma^{tr}_R$ the pseudogroup of local automorphisms of $R\times \mathbb{C}^p\times \mathbb{C}^q$ of the form $F(r,z,w)=(r,f_{l}(r,z,w),f_{tr}(r,w))$, where $f_{tr}$ depending holomorphically on $r$ and $w$ and $f_{l}$ depending holomorphically on $r$ and differentially on $z$ and $w$.\tabularnewline
    \hline
    Family of $h$-deformations & $\Gamma^h_R$ the subpseudogroup of $\Gamma^{tr}_R$ consisting of elements such that $f_{l}$ is also holomorphic on $z$ and $w$. \tabularnewline 
    \hline
    Family of $f$-deformations & $\Gamma^f_R$ the subpseudogroup of $\Gamma_R$ consisting of elements such that $f_{tr}$ do not depend on $r$. \tabularnewline
    \hline
    Family of $f_{supp}$-deformations & $\Gamma^{f_{supp}}_R$ the subpseudogroup of $\Gamma^f_R$ consisting of elements such that $f_{l}$ do not depend on $w$. \tabularnewline
    \hline
 \end{tabular}\\

The coordinates defined from these transition maps are such that $w$ represents the transverse coordinates of the foliations obtained and $z$ coordinates in leaves.

The $f_{supp}$-deformations correspond to the case of a $f$-deformation of a foliation admitting a supplementary foliation. 

\begin{remarque}

The inclusions $\Gamma^{f_{supp}}\subset\Gamma^{f}_R\subset \Gamma^h_R\subset \Gamma^{tr}_R$ show that any family of $f_{supp}$-deformations is a family of $f$-deformations, any family of $f$-deformations is a family of $h$-deformations and any family of $h$-deformations is a family of $tr$-deformations.
\end{remarque}

More concretely, we have :

\begin{definition}
    A family of deformation $(\mathcal{X},R,(X,\mathcal{F}), \gamma_R^{tr})$, denoted $\mathcal{X}/R$ or $\mathcal{X}/(R,\gamma_{R}^{tr})$, of $(X,\F)$ (seen as a transversely holomorphic foliation) parametrized by a complex analytic space $R$, is given by a topological space $\mathcal{X}$, a proper projection $p:\mathcal{X}\to R$ and a $\Gamma^{tr}_R$-structure $\gamma_R^{tr}$ on $\mathcal{X}$ such that, in local coordinates $(r,z,w)$ of this structure, the projection $p$ is the canonical one $p:(r,z,w)\mapsto r$ and the fiber at $0$ of $p$, with the foliation given by $w=cst$, is isomorphic as a foliated complex manifold to $(X,\mathcal{F})$.

    Such a family is called a family of $tr$-deformation of $(X,\F)$ (i.e. deformation of the foliation seen as transversely holomorphic).
\end{definition}

In the following, we will assume that all analytic spaces $R$ we will consider are germs of the form $(R,0)$, and when we will deal with maps between analytic spaces, it should be understood that we are dealing with maps between germs of analytic spaces. 

We also define the notion of a family of $h$-deformations of $(X,\F)$ by considering $\Gamma_R^h$ instead of $\Gamma^{tr}_R$ and that of a family of $f$-deformations of $(X,\F)$ by considering $\Gamma_R^f$ instead of $\Gamma^{tr}_R$. This last type of family corresponds to the situation where the transverse structure is fixed.\\

It is important to note that if we have a family of deformation given by a $\Gamma_R$-structure $\gamma$ on $\mathcal{X}$ and an analytic map $\Theta : R'\to R$, we can define the pull-back $\Theta^* \gamma$ which define a $\Gamma_{R'}$-structure on the fiber product $\mathcal{X}_{\Theta}=\mathcal{X}\times_{R}R'=\{(x,r')\  / \ p(x)=\Theta(r')\}. $\\

We need to know when two families $\mathcal{X}/R$ and $\mathcal{Y}/R'$ are isomorphic families. Without specifying the type of family :

\begin{definition}
    Two families of deformations $\mathcal{X}/R$ and $\mathcal{Y}/R'$ are isomorphic if they define the same $\Gamma$-structure (resp. $\Gamma^f$, $\Gamma^{tr}$), i.e, if there exists $\phi:R\to R'$ isomorphism of analytic spaces and $\Phi:\mathcal{X}\to \mathcal{Y}$ isomorphism of $\Gamma$-structures such that :
    \begin{equation*}
        \begin{tikzcd}
  \mathcal{X} \arrow[rr, "\Phi"] \arrow[rd, ""]
  & & \mathcal{Y} \arrow[ld,"" ]\\
  & R\overset{\phi}{\simeq} R' &  \\ 
\end{tikzcd} 
    \end{equation*}
    commutes.
\end{definition}

\subsection{Multifoliate structures}

The previous $\Gamma$-structures can be described thanks to what Kodaira and Spencer named multifoliate structures in \cite{MUL}. We give here a very brief description of what complex multifoliate structures are.

We denote $(x^{(\alpha)}_1,\dots,x^{(\alpha)}_n)_i$ local coordinates on $M$. Let $(P,\preceq)$ a partially ordered finite set. We say that an integer set is $P$-multifoliated  if it is endowed with a map $\Phi$ :
    \begin{equation*}
        \begin{array}{ccccc}
\Phi & : & \{1,\dots,n\} & \to & P \\
 & & k & \mapsto & [k]:=\Phi(k) \\
\end{array}
    \end{equation*}

    \begin{definition}
    Let $M$ be a differential manifold of dimension $2n$. A multifoliate structure on $M$ correspond to the datum of an atlas $(U_\alpha,(x_1^{(\alpha)},\dots,x_{2n}^{(\alpha)}))_\alpha$, a partially ordered set $(\{\lambda_1,\dots,\lambda_m\}, \succeq)$ and a surjective application $[ \ ]:\{1,\dots,2n\}\to \{\lambda_1,\dots,\lambda_m\}$, such that :
    \begin{equation*}
        \begin{array}{cccc}
            \forall \alpha, \beta, & \frac{\partial x^{(\alpha)}_k}{\partial x^{(\beta)}_l}=0, & if & [k]\cancel{\preceq}[l].
        \end{array}
    \end{equation*}
\end{definition}
    
    When we write for an index $k$ its conjugate $\overline{k}$, we assume that $x^{(\alpha)}_k$ has complex value and that $x^{(\alpha)}_{\overline{k}}:=\overline{x^{(\alpha)}_k}$. Thanks to these “self-conjugated” coordinates $(x^{(\alpha)}_1,\dots,x^{(\alpha)}_n, \overline{x^{(\alpha)}_1},\dots,\overline{x^{(\alpha)}_n})_\alpha$, we can define the notion of holomorphic coordinates, keeping the same formalism by assuming that $[k]\cancel{\preceq}[\overline{l}]$ and $[\overline{k}]\cancel{\preceq}[l]$ for any $k$ and $l$, in other words, non-conjugated integer and conjugated ones are not ordered. Sometimes we will simply denote $(x_1^{(\alpha)},\dots,x_n^{(\alpha)})_\alpha$ what should be thought as "self-conjugated" coordinates.

\begin{exemple}
    \begin{enumerate}
        \item A holomorphic foliation $\F$ on $X$ of codimension $q$ is the datum of an atlas $(U_\alpha,(z^{(\alpha)}_1,\dots,z^{(\alpha)}_{n-q},w^{(\alpha)}_1,\dots,w^{(\alpha)}_q))_\alpha$ such that :
        \begin{equation*}
            \begin{array}{cccccc}
            \forall \alpha, \beta,k,l,k',l' &  \left\{
    \begin{array}{l}
        \frac{\partial w^{(\alpha)}_k}{\partial z^{(\beta)}_{l'}}=\frac{\partial w^{(\alpha)}_k}{\partial \overline{z^{(\beta)}_{l'}}}=\frac{\partial w^{(\alpha)}_k}{\partial \overline{w^{(\beta)}_{k'}}}=0 \\
        \frac{\partial z^{(\alpha)}_k}{\partial \overline{z^{(\beta)}_{l'}}}=\frac{\partial z^{(\alpha)}_k}{\partial \overline{w^{(\beta)}_{k'}}}=0
    \end{array}
\right.  .
        \end{array}
        \end{equation*}
    
    Denoting $z_k^{(\alpha)}=:x^{(\alpha)}_k$, $w_k^{(\alpha)}=:x^{(\alpha)}_{n-q+k}$ and considering $\{\lambda_1,\lambda_2\,\overline{\lambda_1},\overline{\lambda_2}\}$ such that $\lambda_1\preceq \lambda_2$ and $[ \{1,\dots,n-q\}]=\lambda_1$ et $[ \{n-q+1,\dots,n\} ]=\lambda_2$, we get the following multifoliate structure :
    \begin{equation*}
            \begin{array}{cccc}
            \forall \alpha, \beta, & \frac{\partial x^{(\alpha)}_{\lambda}}{\partial x^{(\beta)}_{\lambda'}}=0, & if & [\lambda]\cancel{\preceq}[\lambda'],
        \end{array}
        \end{equation*}
        i.e, 
        \begin{equation*}
            \begin{array}{cccc}
            \forall \alpha, \beta, & \frac{\partial x^{(\alpha)}_{\lambda}}{\partial x^{(\beta)}_{\lambda'}}=0, & if & \left\{
    \begin{array}{l}
        [\lambda]=\lambda_2\mbox{ and }[\lambda']\in\{\lambda_1,\overline{\lambda_1},\overline{\lambda_2}\}\\
        
        [\lambda]=\lambda_1\mbox{ and }[\lambda']\in\{\overline{\lambda_1},\overline{\lambda_2}\}\\

        [\lambda]=\overline{\lambda_2}\mbox{ and }[\lambda']\in\{\lambda_1,\overline{\lambda_1},\lambda_2\}\\

        [\lambda]=\overline{\lambda_1}\mbox{ and }[\lambda']\in\{\lambda_1,\lambda_2\}
    \end{array}
\right. 
        \end{array}
        \end{equation*}
        and in particular,
        \begin{equation*}
            \begin{array}{cccccc}
            \forall \alpha, \beta,k,l,k',l' &  \left\{
    \begin{array}{l}
        \frac{\partial w^{(\alpha)}_k}{\partial z^{(\beta)}_{l'}}=\frac{\partial w^{(\alpha)}_k}{\partial \overline{z^{(\beta)}_{l'}}}=\frac{\partial w^{(\alpha)}_k}{\partial \overline{w^{(\beta)}_{k'}}}=0 \\
        \frac{\partial z^{(\alpha)}_k}{\partial \overline{z^{(\beta)}_{l'}}}=\frac{\partial z^{(\alpha)}_k}{\partial \overline{w^{(\beta)}_{k'}}}=0
    \end{array}
\right.  .
        \end{array}
        \end{equation*}
        \item The datum of two holomorphic foliations $\F$ and $\mathcal{G}$ everywhere transverse and supplementary ($TX=T\F\oplus T\mathcal{G}$) on $X$, corresponds to the datum of a particular atlas $(U_\alpha,(z^{(\alpha)}_1,\dots,z^{(\alpha)}_{p},w^{(\alpha)}_1,\dots,w^{(\alpha)}_q))_\alpha$ such that :
          \begin{equation*}
            \begin{array}{cccccc}
            \forall \alpha, \beta,k,l,k',l' &  \left\{
    \begin{array}{l}
        \frac{\partial w^{(\alpha)}_k}{\partial z^{(\beta)}_{l'}}=\frac{\partial w^{(\alpha)}_k}{\partial \overline{z^{(\beta)}_{l'}}}=\frac{\partial w^{(\alpha)}_k}{\partial \overline{w^{(\beta)}_{k'}}}=0 \\
        \frac{\partial z^{(\alpha)}_l}{\partial \overline{z^{(\beta)}_{l'}}}=\frac{\partial z^{(\alpha)}_l}{\partial \overline{w^{(\beta)}_{k'}}}=\frac
        {\partial z^{(\alpha)}_l}{\partial w^{(\beta)}_{k'}}=0
    \end{array}
\right.  .
        \end{array}
        \end{equation*}
        Denoting $x^{(\alpha)}_k:=z^{(\alpha)}_k$, $x^{(\alpha)}_{p+k}:=w^{(\alpha)}_k$ and considerating $\{\lambda_1,\lambda_2,\overline{\lambda_1},\overline{\lambda_2}\}$ without ordering any element to another and given $[\{1,\dots,p\}]=\lambda_1$ and $[\{p+1,\dots,n\}]=\lambda_2$, we get the following multifoliate structure :
        \begin{equation*}
            \begin{array}{cccc}
            \forall \alpha, \beta, & \frac{\partial x^{(\alpha)}_{k}}{\partial x^{(\beta)}_l}=0, & if & [k]\cancel{\preceq}[l],
        \end{array}
        \end{equation*}
        i.e,
                \begin{equation*}
            \begin{array}{cccc}
            \forall \alpha, \beta, & \frac{\partial x^{(\alpha)}_{\lambda}}{\partial x^{(\beta)}_{\lambda'}}=0, & if & \left\{
    \begin{array}{l}
        [\lambda]=\lambda_2\mbox{ and }[\lambda']\in\{\lambda_1,\overline{\lambda_1},\overline{\lambda_2}\}\\
        
        [\lambda]=\lambda_1\mbox{ and }[\lambda']\in\{\overline{\lambda_1},\overline{\lambda_2},\lambda_2\}\\

        [\lambda]=\overline{\lambda_2}\mbox{ and }[\lambda']\in\{\lambda_1,\overline{\lambda_1},\lambda_2\,\}\\

        [\lambda]=\overline{\lambda_1}\mbox{ and }[\lambda']\in\{\lambda_1,\lambda_2,\overline{\lambda_2}\}
    \end{array}
\right. 
        \end{array}
        \end{equation*}
        and in particular,
        \begin{equation*}
            \begin{array}{cccccc}
            \forall \alpha, \beta,k,l,k',l' &  \left\{
    \begin{array}{l}
        \frac{\partial w^{(\alpha)}_k}{\partial z^{(\beta)}_{l'}}=\frac{\partial w^{(\alpha)}_k}{\partial \overline{z^{(\beta)}_{l'}}}=\frac{\partial w^{(\alpha)}_k}{\partial \overline{w^{(\beta)}_{k'}}}=0 \\
        \frac{\partial z^{(\alpha)}_l}{\partial \overline{z^{(\beta)}_{l'}}}=\frac{\partial z^{(\alpha)}_l}{\partial \overline{w^{(\beta)}_{k'}}}=\frac
        {\partial z^{(\alpha)}_l}{\partial w^{(\beta)}_{k'}}=0
    \end{array}
\right.  .
        \end{array}
        \end{equation*}
        \end{enumerate}
\end{exemple}
    
\subsection{Almost multifoliate structures}

Another way of defining multifoliate structures is by means of differential complex subundles of the complexified tangent bundle $T_\C M:=TM\otimes \C$ of $M$. Since the aim is to define local systems of differential equations, one natural thing to do is to consider the subundles whose fields define the equations of interest. We then have :

\begin{definition}
    Let $M$ be a differential manifold (compact). We call almost complex multifoliate structure the datum of a collection $(F_{\iota_1},\dots,F_{\iota_m})$ of subundles of $T_\C M$ and a partial order relationship $\leq$ on $\{\iota_1,\dots,\iota_m\}$, such that :
    \begin{enumerate}[label=(\roman*)]
        \item For all $k\in\{1,\dots,m\}$, $T_\C M=F_{\iota_k}+\overline{F_{\iota_k}}$
        \item For all $k,l\in\{1,\dots,m\}$, $F_{\iota_k}\subset F_{\iota_l}$, if and only if, $\iota_k\leq \iota_l$.
    \end{enumerate}
\end{definition}

\begin{remarque}
    The partial order relationship of this definition does not correspond to that of the definition given in terms of local coordinates.
\end{remarque}

\begin{exemples}\label{exemplealmost}
    \begin{enumerate}
        \item Any almost complex structure $T_\C M=T^{0,1}\oplus T^{1,0}$ on $M$.
        \item Any transversaly holomorphic foliation on $M$ can be defined thanks to the datum of a unique subundle $F\subset T_\C M$ defined as $ker(T_\C M\xrightarrow[]{pr_{\nu\F^{1,0}}}\nu\F^{1,0})$. It is such that $F+\overline{F}=T_\C M$. For more details see (\cite{NIC}, 1).
        \item Any holomorphic foliation on $M$ is defined thanks to the datum of a unique pair $(T^{0,1},F)$ defining respectively the complex structure on $M$ and the transversaly holomorphic structure of the foliation. The guarantee that there are coordinates in the leaves adapted to the global complex structure comes from the relationship $T^{0,1}\subset F$.
        \item Two supplementary holomorphic foliations on $M$ are defined thanks to the datum of $(T^{0,1},F,G)$ defining respectively the complex structure on $M$ and the transversaly holomorphic structure of each foliation. The guarantee that the coordinates are adapted to the global complex structure comes from the relationship $T^{0,1}\subset F$ and $T^{0,1}\subset G$.
    \end{enumerate}
\end{exemples}

\begin{definition}
   An almost complex multifoliate structure $(F_{\iota_1},\dots,F_{\iota_m})$ is integrable when each $F_{\iota_k}$ is involutive.  
\end{definition}

We can show that any multifoliate structure on $M$ comes from a unique almost complex multifoliate structure integrable. It is a part of the theorem \ref{representationfamille}.\\

\subsection{Deformation Theory}\label{deformationtheory}

Until now, we have seen three different ways to define our structures of interest : holomorphic and transversely holomorphic foliations.\\

\underline{As $\Gamma$-structures :} In the table given in \ref{Gammastructures} we already described what a family of deformations is from the $\Gamma$-structure point of view for different type of deformations.\\

\underline{As multifoliate structures :} From the point of view of multifoliate structures, a family of deformations parametrized on an analytic space $R$ corresponds to the data of a family of atlases $(U_\alpha,(x_1^{(\alpha)}(r),\dots,x_{2n}^{(\alpha)}(r))_{r\in R}$, a partially ordered set $(\{\lambda_1,\dots,\lambda_m\}, \succeq)$ and a surjective application $[ \ ]:\{1,\dots,2n\}\to \{\lambda_1,\dots,\lambda_m\}$, such that the coordinates are holomorphic as functions of $r$ and  :
    \begin{equation*}
        \begin{array}{cccc}
            \forall r,\alpha, \beta, & \frac{\partial x^{(\alpha)}_k(r)}{\partial x^{(\beta)}_l(r)}=0, & if & [k]\cancel{\preceq}[l].
        \end{array}
    \end{equation*}\\

\underline{As integrable almost complex multifoliate structures :} From the point of view of integrable almost complex structures, a family of deformations parametrized on an analytic space $R$ corresponds to the data of a partial order relationship $\leq$ on $\{\iota_1,\dots,\iota_m\}$ and a collection $(F_{\iota_1}(r),\dots,F_{\iota_m}(r))_{r\in R}$ of subundles of $T_\C M$ obtained thanks to a family of endomorphisms holomorphic in $r\in R$ :
\begin{equation*}
    \Phi:R\times T_\C M\to T_\C M,
\end{equation*}
as $F_{\iota_k}(r)=\Phi_r(F_{\iota_k})$, where $\Phi_r:=\Phi(r,\cdot)$, and such that :

\begin{equation*}
    \forall r, \ F_{\iota_k}(r)\subset F_{\iota_l}(r) \iff \iota_k\leq \iota_l.
\end{equation*}
and all almost complex multifoliate structures defined are integrable. Since we deal with complex structures, we assume that $\Phi_r|_{T^{0,1}}=\overline{\Phi_r|_{T^{1,0}}}$.\\

Since we gave different type of deformations for holomorphic foliations, we will recap the different ones for each point of view in the following table. Just before, let us note two things for a holomorphic foliation given as $(T^{0,1},F)$.\\

Firstly, according to the table of \ref{Gammastructures}, a $f$-deformation $(T^{0,1}_r,F_r)_{r\in R}$ is such that the transverse structure is the same. It means that $F_r=F_0=F$. We need also that $T^{0,1}_r\subset F$. In fact, we have :
\begin{lemme}
    The bundle $T^{0,1}_r$ is a subbundle of $F$ if and only if $(\Phi_r-id)|_{T^{0,1}}$ takes values in $T\mathcal{F}^{1,0}$.
\end{lemme}

\begin{preuve}
    Thanks to $T\mathcal{F}^{1,0}=T^{1,0}\cap F$, the converse is immediate.

    When $T^{0,1}_r\subset F$, the direct implication comes from the fact that $(\Phi_r-id)|_{T^{0,1}}$ has to takes values in $F$ since $T^{0,1}\subset F$, and then take values in $T^{1,0}\cap F=T\mathcal{F}^{1,0}$.
\end{preuve}

Secondly, the datum of two supplementary holomorphic foliations $\F, \mathcal{G}$ on $X$ (i.e $T_\C M=T\F^{0,1}\oplus T\F^{1,0}\oplus T\mathcal{G}^{0,1}\oplus T\mathcal{G}^{1,0}$) is such that $\F$ is given by $(T^{0,1},F)$ and $\mathcal{G}$ by $(T^{0,1},G)$. But we can say more thanks to the definition of $F$ and $G$ : $F=T\F^{1,0}\oplus T\F^{0,1}\oplus T\mathcal{G}^{0,1}$ and $G=T\mathcal{G}^{1,0}\oplus T\mathcal{G}^{0,1}\oplus T\F^{0,1}$. A $f_{supp}$-deformation of $\F$ then corresponds only to a deformation of $T\F^{0,1}$ (and $T\F^{0,1}$), obtained with $(\Phi_r-id)|_{T^{0,1}}:T_\C M\to T_\F^{1,0}$ and $(\Phi_r-id)|_{T\mathcal{G}^{0,1}}\equiv 0$.

\begin{tabular}{|p{2cm}|p{4cm}|p{4cm}|p{4cm}|}

    \hline
    \textbf{Family of} & \textbf{$\Gamma$-structures} & \textbf{Multifoliate structures} & \textbf{Almost Complex multifoliate structures} \tabularnewline
    \hline
    $tr$-deformations (the foliation is seen as a transversaly holomorphic one) & $\Gamma^{tr}_R$ the pseudogroup of local automorphisms of $ R\times \mathbb{R}^{2n-2q}\times \mathbb{C}^q$ of the form $F(r,x,y,w)=(r,f_{l}(r,x,y,x',y'),$ $f_{tr}(r,w))$, where $f_{tr}$ depending holomorphically on $r$ and $w:=x'+iy'$ and $f_{l}$ depending holomorphically on $r$ and differntially on $x,y,x'$ and $y'$. & Atlas of the form $(x_l^{(\alpha)}(r), w_k^{(\alpha)}(r))$ with $(x^{(\alpha)}_l(r))_{1\leq l\leq 2n-2q} \in \R^{2n-2q}$ are the differential coordinates on the leaves and $(w^{(\alpha)}_k)_{1\leq k \leq q}\in \C^q $ such that  $\frac{\partial w^{(\alpha)}_k(r)}{\partial x^{(\beta)}_l(r)}=0$. and $\frac{\partial w^{(\alpha)}_k(r)}{\partial \overline{w^{(\beta)}_l(r)}}=0$ & $(F_r)$ with $F_r$ involutives and $F_r+\overline{F_r}=T_\C M$. \tabularnewline
    \hline
    $h$-deformations & $\Gamma^h_R$ the subpseudogroup of $\Gamma^{tr}_R$ consisting of elements such that $f_{l}$ is also holomorphic on $z$ and $w$. & Atlas of the form $(z_l^{(\alpha)}(r), w_k^{(\alpha)}(r))$ with $(z^{(\alpha)}_l(r))_{1\leq l\leq p} \in \C^{p}$ and $(w^{(\alpha)}_k)_{1\leq k \leq q}\in \C^q $ such that the coordinates are holomorphic and $\frac{\partial w^{(\alpha)}_k(r)}{\partial z^{(\beta)}_l(r)}=0$.& $(T^{0,1}_r,F_r)$, with $T^{0,1}_r$ and $F_r$ involutives, $F_r+\overline{F_r}=T^{0,1}_r\oplus\overline{T^{0,1}_r}=T_\C M.$ \tabularnewline 
    \hline
    $f$-deformations & $\Gamma^f_R$ the subpseudogroup of $\Gamma_R$ consisting of elements such that $f_{tr}$ do not depend on $r$ & Atlas of the form $(z_l^{(\alpha)}(r), w_k^{(\alpha)}(r))$ with $(z^{(\alpha)}_l(r))_{1\leq l\leq p} \in \C^{p}$ and $(w^{(\alpha)}_k)_{1\leq k \leq q}\in \C^q $ such that the coordinates are holomorphic and $\frac{\partial w^{(\alpha)}_k(r)}{\partial z^{(\beta)}_l(r)}=0$. & $(T^{0,1}_r,F)$, with $T^{0,1}_r$ and $F$ involutives and $F+\overline{F}=T^{0,1}_r\oplus\overline{T^{0,1}_r}=T_\C M$, obtained with $\Phi_r:T_\C M\to T_\C M$ such that $(\Phi_r-id)|_{T^{0,1}}:T^{0,1}\to T\F^{1,0}$.  \tabularnewline
    \hline
    $f_{supp}$-deformations & $\Gamma^{f_{supp}}_R$ the subpseudogroup of $\Gamma^f_R$ consisting of elements such that $f_{l}$ do not depend on $w$. & Atlas of the form $(z_l^{(\alpha)}(r), w_k^{(\alpha)}(r))$ with $(z^{(\alpha)}_l(r))_{1\leq l\leq p} \in \C^{p}$ and $(w^{(\alpha)}_k)_{1\leq k \leq q}\in \C^q $ such that the coordinates are holomorphic and $\frac{\partial w^{(\alpha)}_k(r)}{\partial z^{(\beta)}_l(r)}=0$. & $(T^{0,1}_r,F)$ with $T^{0,1}_r=T\F^{0,1}_r\oplus T\mathcal{G}^{0,1}$, $T^{0,1}_r$ and $F$ involutives and $F+\overline{F}=T^{0,1}_r\oplus\overline{T^{0,1}_r}=T_\C M$, obtained with $\Phi_r:T_\C M\to T_\C M$ such that $(\Phi_r-id)|_{T^{0,1}}:T\F^{0,1}\to T\F^{1,0}$ and $(\Phi_r-id)|_{T\mathcal{G}^{0,1}}\equiv0$.\tabularnewline
    \hline
 \end{tabular}\\

 Moreover, a family of deformations of multifoliate structures can be represented by a global derivations (i.e a linear map $u$ between the spaces of global differential forms such that $u(A^*_M)\subset A^{*+s}_M$ for a certain $s$ and verifying Leibniz formula) itself represented by what Kodaira and Spencer named in \cite{MUL} (I.1) jet forms. We describe briefly what jet forms are, it will be useful later. For any details about the relation between derivations and jet forms or the definition of the bracket for jet forms, see \cite{MUL}.

 \begin{definition}
A jet form of degree $r$ (or $r$-jet form) on $M$ is the data of a covering $(U_\alpha)_\alpha$ of $M$ and a collection $\{(\phi^{(\alpha)},\xi^{(\alpha)})\}_{\alpha\in I}$ of pairs $(\phi^{(\alpha)},\xi^{(\alpha)})$ of vector forms $\phi^{(\alpha)}$ and $\xi^{(\alpha)}$ of respective degrees $r$ and $r+1$, defined on $U_\alpha$ such that, on intersections $U_\alpha\cap U_\beta$ we have :
\begin{equation*}
    \left\{
    \begin{array}{ll}
        \phi^{(\beta)}=\phi^{(\alpha)}\\
        \xi^{(\beta)}-d\phi^{(\beta)}=\xi^{(\alpha)}-d\phi^{(\alpha)}
    \end{array}.
    \right.
\end{equation*}
\end{definition}

\begin{remarque}  
The quantities $d\phi^{(\beta)}$ and $d\phi^{(\alpha)}$ are not necessarily equal, they depend on the local charts. In particular, the first equality does not necessarily imply $d\phi^{(\beta)}=d\phi^{(\alpha)}$.
\end{remarque}

Finally, a jet form is said to respect a given multifoliate structure when :

\begin{definition}
Let $u:=\{(\phi^{(\alpha)},\xi^{(\alpha)})\}_\alpha$ be a $p$-jet form on $M$. We say that $u$ respects the $\Gamma_P$-structure of $M$ when :
\begin{equation*}
    \xi^{(\alpha)}_{\lambda,\lambda_1,\dots,\lambda_{p+1}}=0, \ \mbox{for } \lambda_1\cancel{\succeq} \lambda \mbox{ and } \dots \mbox{ and } \lambda_{p+1}\cancel{\succeq} \lambda,
\end{equation*}

where $\xi^{(\alpha)}=\sum\limits_{\lambda,\lambda_1,\dots,\lambda_{p+1}} \xi^{(\alpha)}_{\lambda,\lambda_1,\dots,\lambda_{p+1}}dx^{(\alpha)}_{\lambda_1}\wedge\dots \wedge dx^{(\alpha)}_{\lambda_{p+1}} \otimes \frac{\partial}{\partial x^{(\alpha)}_\lambda}$.
\end{definition}

Kodaira and Spencer showed the following (note that $d$ as a derivation can be represented by a jet form) :

\begin{theoreme}(Kodaira \& Spencer)\label{representationfamille}
Any family of deformations of multifoliate structures $(U_\alpha,(x_1^{(\alpha)}(r),\dots,x_{2n}^{(\alpha)}(r))_{r\in R}$, $(\{\lambda_1,\dots,\lambda_m\}, \succeq)$, determines a family of global $1$-jet form $\{v(r)=(\Phi(r),\Xi(r)), \ r\in R\}$ respecting the multifoliate structure and verifying :
\begin{equation*}
    \left\{
    \begin{array}{ll}
        [v,v]=0\\
        v(0)=d
    \end{array},
    \right.
\end{equation*}
and for $\Phi(r)=\sum\limits_\lambda \psi^{(\alpha)}_\lambda (x,r) dx^{(\alpha)}_\lambda(0)$, we have :
\begin{equation*}
    \psi^{(\alpha)}_{\lambda_1}(r)\cdot  x^{(\alpha)}_{\lambda_2}(r) =0, \ \mbox{ for } \lambda_1\cancel{\succeq} \lambda_2.
\end{equation*}
Reciprocally, any such 1-jet form determines a family of deformations of the initial structure $(U_\alpha,(x_1^{(\alpha)}(0),\dots,x_{2n}^{(\alpha)}(0))$, $(\{\lambda_1,\dots,\lambda_m\}, \succeq)$.
\end{theoreme}

In the theorem, the global $1$-form vector $\Phi(r)$ is exactly the family of endomorphisms $\Phi(r):T_\C M\to T_\C M$. The second term of the pair $\Xi(r)$ allows to store the involutivity conditions when the condition $[v,v]=0$ is verified.  

\begin{remarque}
    Considering the representation of families of small deformations by means of a description in terms of jet forms makes it possible, when it exists, to make natural use of the existence of a foliation $\mathcal{G}$ everywhere transverse and supplementary to $\mathcal{F}$.
\end{remarque}

We give some ideas that will enlighten how are related the description of a family of deformations, as given by suitable local coordinates, to the description obtained by means of $1$-jets form. We will call the transition from one to the other the Kodaira-Spencer procedure. We will refer to this procedure later.\\

Since the general procedure is the same, let us assume, for example, that a set : 
\begin{equation*}
(U_\alpha,(x^{(\alpha)}_1(r),\dots,x^{(\alpha)}_{n} (r)))_{\alpha}
\end{equation*}

of local coordinates adapted to a family of deformations of a complex multifoliate structure is given. We denote $x^{(\alpha)}_k$ for $x^{(\alpha)}_k(0)$. So we are looking to construct a global $1$-jet form that respects the multifoliate structure. The challenge is to build a map depending on the parameter $r$ that contains the way in which a base of fields is transformed.

First, we denote $(g^{(\alpha,\beta)}_{k,l}(r))_{\alpha,\beta,k,l}$ the transition maps :

\begin{equation*}
    dx^{(\alpha)}_k(r)=\sum\limits_q g^{(\alpha,\beta)}_{k,l}(r) dx^{(\beta)}_l(r),
\end{equation*}

which are such that :
\begin{equation*}
    g^{(\alpha,\beta)}_{k,l}(r)=0 \mbox{ for }l\cancel{\succeq} k,
\end{equation*} 
by definition of multifoliate structures.

Denoting then $1=\sum\limits_\beta \rho_\beta (x)$ a partition of unity, which we assume without loss of generality is adapted to $(U_\alpha)_\alpha$ (in particular, these functions are independent of $r$ and therefore holomorphic with respect to this variable), we set :

\begin{equation*}
    b^{(\alpha)}_{k,l}(r)=\sum\limits_\beta \rho_\beta\sum\limits_{p} g^{(\alpha,\beta)}_{k,p}(r)g^{(\beta,\alpha)}_{p,l}(0),
\end{equation*}

and :

\begin{equation*}
    \omega^{(\alpha)}_k(r):=\sum\limits_l b^{(\alpha)}_{l,k}(r)\frac{\partial }{\partial x^{(\alpha)}_l(r)}.
\end{equation*}

and finally :

\begin{equation*}
    \omega^{(\alpha)}(r):=\sum\limits_k \psi^{(\alpha)}_k(r) dx^{(\alpha)}_k.
\end{equation*}

We can show (\cite{MUL}, 8 and 5) that the $(\omega^{(\alpha)}+\overline{\omega^{(\alpha)}})_\alpha$ can be glued to obtain a global section $\psi$ and that the $(\omega^{(\alpha)}_k,\overline{\omega_k^{(\alpha)}})_k=(\omega^{(\alpha)}_\lambda)_{\lambda=k,\overline{k}}$ are forming local basis of $T_\C M$. Note that $\psi$ can be seen as a map which for a vector field of the initial structure gives a vector field dependent on the parameter $r$ (and coinciding with the initial vector field at $0$). The vector fields obtained correspond to a basis of the tangent bundle of the structures defined in $r$ and are respecting the multifoliate structure.

So we define a global vector $1$-form. We still have to define a global $2$-form verifying some properties to obtain a $1$-jet form that respects the multifoliate structure and representing the starting family. 

From what we said above, we can express the Poisson brackets of the vector fields $\psi^{(\alpha)}_{\lambda}$ and $\psi^{(\alpha)}_{\lambda '}$ as :

\begin{equation*}
    [\psi^{(\alpha)}_\lambda,\psi^{(\alpha)}_{\lambda '}]=:\sum\limits_{\mu} c^{(\alpha)}_{\lambda,\lambda ',\mu}(r) \psi^{(\alpha)}_\mu.
\end{equation*}

Then, by definition of the bracket and of multifoliate structures by means of pseudogroups :

\begin{equation*}
    c^{(\alpha)}_{\lambda,\lambda ',\mu}(r)=0, \mbox{ for } \lambda\cancel{\succeq} \mu \mbox{ and } \lambda '\cancel{\succeq} \mu, 
\end{equation*}

and we define on each $U_\alpha$, a local vector 2-form :

\begin{equation*}
    c_\alpha(r):=\sum\limits_{\mu} (\sum\limits_{\lambda,\lambda '} c^{(\alpha)}_{\lambda,\lambda ',\mu} dx^{(\alpha)}_\lambda\wedge dx^{(\alpha)}_{\lambda'})\frac{\partial }{\partial x^{(\alpha)}_{\mu}}.
\end{equation*}

We can then check that $u(r):=(\psi_\alpha,c_\alpha)_\alpha$ defines a global $1$-jet form respecting the multifoliate structure and such that $[u(r),u(r)]=0$. This last equation corresponds to an integrability condition of the Frobenius type. By calculating this bracket, we can show that, under this hypothesis, bundles ($T^{0,1}$, $F$, ... depending on the multifoliate structure) are indeed involutives. 

Reciprocally, by virtue of the latter integrability condition and a Frobenius-type theorem, we can show that from such a jet form there exists local coordinates $(U_\alpha,(x^{(\alpha)}_1(r),\dots,x^{(\alpha)}_n(r)))_\alpha$ solutions of :

\begin{equation*}
    \psi^{(\alpha)}_{\lambda '}(r) \cdot x^{(\alpha)}_\lambda(r)=0, \mbox{ for } \lambda '\cancel{\succeq} \lambda.
\end{equation*}

\subsection{Sheaves of infinitesimal deformations}

We can define an application representing, for a given family of deformations, the infinitesimal deformation obtained in a given direction of the parameter space. Such a map is called a Kodaira-Spencer map. Once we identify what infinitesimal deformations are for a given structure, the aim is to construct a family of deformation that, in a sense, integrates this infinitesimal deformation.

With atlases description, we can define such a map in the case of a family of $tr$-deformations $\mathcal{X}/R$ locally given by $(U_\alpha,(r,x^{(\alpha)}_1,\dots,x_{2n-2q}^{(\alpha)},w^{(\alpha)}_1,\dots,w^{(\alpha)}_q))_\alpha$ in the following way. 

The transition maps on $U^{(\alpha,\beta)}$ are given by :
\begin{equation*}
    F^{(\alpha,\beta)}(r,x^{(\alpha)},w^{(\alpha)})=(r,f^{(\alpha,\beta)}_\mathcal{L}(r,x^{(\alpha)},w^{(\alpha)}),f^{(\alpha,\beta)}_{tr}(r,w^{(\alpha)})).
\end{equation*}
To account, in a given direction, for the deformation obtained with the family $\mathcal{X}/R$, we calculate for any $\frac{\partial}{\partial r}\in T_0 R$, the tangent Zariski space of $R$ at $0$, the quantities :

\begin{equation*}
    \widehat{\theta_{\alpha,\beta}}:=\sum\limits_{k=1}^p \frac{\partial f^{(\alpha,\beta)}_{\mathcal{L},k}}{\partial r}|_{r=0}\cdot \frac{\partial }{\partial x^{(\beta)}_k}+\sum\limits_{k=1}^q \frac{\partial f^{(\alpha,\beta)}_{tr,k}}{\partial r}|_{r=0}\cdot \frac{\partial }{\partial w^{(\beta)}_k}.
\end{equation*}

In particular, these fields are $\mathcal{C}^\infty$ and their flows define local isomorphisms of $\mathcal{F}$ as a transversaly holomorphic foliation. We denote $\widehat{\Theta_{tr}}$ the sheaf of germs of such fields. 

From the above equalities, we obtain :

\begin{equation*}
    f^{(\alpha,\beta)}_{\mathcal{L},k}(r,x^{(\alpha)},w^{(\alpha)})=f^{(\gamma,\beta)}_{\mathcal{L},k}(F^{(\alpha,\gamma)}(r,x^{(\alpha)},w^{(\alpha)}))
\end{equation*}

and :

\begin{equation*}
    f^{(\alpha,\beta)}_{tr,k}(r,w^{(\alpha)})=f^{(\gamma,\beta)}_{tr,k}(r,f^{(\alpha,\gamma)}_{tr,1}(r,w^{(\alpha)}),\dots,f^{(\alpha,\gamma)}_{tr,q}(r,w^{(\alpha)}))
\end{equation*}

we respectively obtain :

\begin{equation*}
    \frac{\partial f^{(\alpha,\beta)}_{\mathcal{L},k}}{\partial r}=\frac{\partial f^{(\gamma,\beta)}_{\mathcal{L},k}}{\partial r}+\sum\limits_{l=1}^p \frac{\partial f^{(\gamma,\beta)}_{\mathcal{L},k}}{\partial x^{(\gamma)}_l}\frac{\partial f^{(\alpha,\gamma)}_{\mathcal{L},l}}{\partial r}+\sum\limits_{l=1}^q \frac{\partial f^{(\gamma,\beta)}_{\mathcal{L},k}}{\partial w^{(\gamma)}_l}\frac{\partial f^{(\alpha,\gamma)}_{tr,l}}{\partial r}
\end{equation*}

and :

\begin{equation*}
     \frac{\partial f^{(\alpha,\beta)}_{tr,k}}{\partial r}=\frac{\partial f^{(\gamma,\beta)}_{tr,k}}{\partial r}+ \sum\limits_{l=1}^q \frac{\partial f^{(\gamma,\beta)}_{tr,k}}{\partial w^{(\gamma)}_l} \frac{\partial 
 f^{(\alpha,\gamma)}_{tr,l}}{\partial r}
\end{equation*}

and therefore, respectively : 

\begin{equation*}
\begin{array}{cccc}
     \sum\limits_{k=1}^p \frac{\partial f^{(\alpha,\beta)}_{\mathcal{L},k}}{\partial r}|_{r=0}\cdot \frac{\partial }{\partial x^{(\beta)}_k} & = & & \sum\limits_{k=1}^p \frac{\partial f^{(\gamma,\beta)}_{\mathcal{L},k}}{\partial r}|_{r=0}\cdot \frac{\partial }{\partial x^{(\beta)}_k}\\
     & & + &  \sum\limits_{k=1}^p\sum\limits_{l=1}^p \frac{\partial f^{(\gamma,\beta)}_{\mathcal{L},k}}{\partial x^{(\gamma)}_l}\frac{\partial f^{(\alpha,\gamma)}_{\mathcal{L},l}}{\partial r}|_{r=0}\cdot \frac{\partial }{\partial x^{(\beta)}_k} \\
     & & + & \sum\limits_{k=1}^p\sum\limits_{l=1}^q \frac{\partial f^{(\gamma,\beta)}_{\mathcal{L},k}}{\partial w^{(\gamma)}_l}\frac{\partial f^{(\alpha,\gamma)}_{tr,l}}{\partial r}|_{r=0}\cdot \frac{\partial }{\partial x^{(\beta)}_k}\\
     
    & = & & \sum\limits_{k=1}^p \frac{\partial f^{(\gamma,\beta)}_{\mathcal{L},k}}{\partial r}|_{r=0}\cdot \frac{\partial }{\partial x^{(\beta)}_k}\\
     & & + &  \sum\limits_{l=1}^p \frac{\partial f^{(\alpha,\gamma)}_{\mathcal{L},l}}{\partial r}|_{r=0}\cdot \frac{\partial }{\partial x^{(\gamma)}_l} \\
     & & + & \sum\limits_{k=1}^p\sum\limits_{l=1}^q \frac{\partial f^{(\gamma,\beta)}_{\mathcal{L},k}}{\partial w^{(\gamma)}_l}\frac{\partial f^{(\alpha,\gamma)}_{tr,l}}{\partial r}|_{r=0}\cdot \frac{\partial }{\partial x^{(\beta)}_k}
\end{array}
\end{equation*}

and :

\begin{equation*}
\begin{array}{cccc}
      \sum\limits_{k=1}^q \frac{\partial f^{(\alpha,\beta)}_{tr,k}}{\partial r}|_{r=0}\cdot \frac{\partial }{\partial w^{(\beta)}_k} & = & & \sum\limits_{k=1}^q \frac{\partial f^{(\gamma,\beta)}_{tr,k}}{\partial r}|_{r=0}\cdot \frac{\partial }{\partial w^{(\beta)}_k} \\
      &  & + & \sum\limits_{k=1}^q \sum\limits_{l=1}^q \frac{\partial f^{(\gamma,\beta)}_{tr,k}}{\partial w^{(\gamma)}_l} \frac{\partial f^{(\alpha,\gamma)}_{tr,l}}{\partial r}|_{r=0}\cdot \frac{\partial }{\partial w^{(\beta)}_k}\\
      & = &  & \sum\limits_{k=1}^q \frac{\partial f^{(\gamma,\beta)}_{tr,k}}{\partial r}|_{r=0}\cdot \frac{\partial }{\partial w^{(\beta)}_k}\\
      & & + & \sum\limits_{k=1}^q \frac{\partial f^{(\alpha,\gamma)}_{tr,k}}{\partial r}|_{r=0}\cdot \frac{\partial }{\partial w^{(\gamma)}_k}\\
      & & - & \sum\limits_{k=1}^p \sum\limits_{l=1}^q \frac{\partial f^{(\gamma,\beta)}_{\mathcal{L},k}}{\partial w^{(\gamma)}_l} \frac{\partial f^{(\alpha,\gamma)}_{tr,l}}{\partial r}|_{r=0}\cdot \frac{\partial }{\partial x^{(\beta)}_k}.
\end{array}
\end{equation*}

We then note that :

\begin{equation*}
    \widehat{\theta_{\alpha,\beta}}=\widehat{\theta_{\gamma,\beta}}+\widehat{\theta_{\alpha,\gamma}}.
\end{equation*}

 In other words, $(\widehat{\theta_{\alpha,\beta}})_{\alpha,\beta}$ is a 1-cocycle of the sheaf $\widehat{\Theta_{tr}}$. In particular, it defines a class in $H^1(M,\widehat{\Theta_{tr}})$. Finaly, we define a Kodaira-Spencer map :

\begin{equation*}
    \rho_{tr}: T_0 R\to H^1(M,\widehat{\Theta_{tr}})
\end{equation*}

 However, we may consider that the leaf structure is of no interest to us during transverse holomorphic deformation, and consider only the infinitesimal deformation of the complex structure obtained thanks to the submersions defining transverse holomorphic foliation. To do this, we consider only the deformation of complex coordinates $(w^{(\alpha)})$ by calculating for $\frac{\partial}{\partial r}\in T_0 R$, the quantities :

 \begin{equation*}
    \theta_{\alpha,\beta}:=\sum\limits_{k=1}^q \frac{\partial f^{(\alpha,\beta)}_{tr,k}}{\partial r}|_{r=0}\cdot \frac{\partial }{\partial w^{(\beta)}_k}.
\end{equation*}

The same considerations lead us to the definition of a second Kodaira-Spencer map. Denoting $\Theta_{tr}$ the quotient of $\widehat{\Theta_{tr}}$ by the sheaf $\mathcal{A}^0(T\F)$ of germs of differential vector fields tangent to the leaves, we get :

\begin{equation*}
    0 \xrightarrow[]{} \mathcal{A}^0(T\F) \xrightarrow[]{} \widehat{\Theta_{tr}} \xrightarrow[]{} \Theta_{tr} \xrightarrow[]{} 0,
\end{equation*}

and we can define :

\begin{equation*}
    \rho_{tr}': T_0 R\to H^1(M,\Theta_{tr})
\end{equation*}

It is worth noting at this point that giving two different maps that a priori account for the same infinitesimal deformations is not problematic. In fact, since $\Theta_{tr}$ is the quotient of $\widehat{\Theta_{tr}}$ by an acyclic subsheaf ($\mathcal{C}^\infty$-module), we get for $k>0$ :

\begin{equation*}
    H^k(M,\Theta_{tr})\simeq H^k(M,\widehat{\Theta_{tr}}).
\end{equation*}

We can do the same in the case of holomorphic foliations and define a Kodaira-Spencer map associated with a holomorphic family of deformations of the foliation :

\begin{equation*}
    \rho_{h}:T_0 R\to H^1(M,\Theta_h),
\end{equation*}

where $\Theta_h$ is the sheaf of germs of vector fields whose flows are local isomorphisms of the holomorphic foliation (in particular, they care about preserving the holomorphic character of the leaves).

It is also natural to consider the subsheaf of germs of holomorphic vector fields tangent to the leaves $T\F^{1,0}$. The link between $T\F$, $\Theta_h$ and $\Theta_{tr}$ is given by :

\begin{equation*}
    0\xrightarrow[]{} T\F^{1,0} \xrightarrow[]{} \Theta_h\xrightarrow[]{} \Theta_{tr}\xrightarrow[]{} 0.
\end{equation*}

This is summarized in the following table :\\

\begin{tabular}{|p{5cm}|p{5cm}|p{5cm}|}
    \hline
    \textbf{Structures} & \textbf{Sheaves associated to the infinitesimal deformations} & \textbf{Kodaira-Spencer map}\tabularnewline
    \hline
    
    Family of $tr$-deformations & $\Theta_{tr}$ locally given by classes of vector fields of the form :
    \begin{equation*}
        [\sum\limits_k b_{tr, k}(r,w^{(\alpha)})\frac{\partial}{\partial w^{(\alpha)}_k}],
    \end{equation*}
    with $b_{tr, k}$ holomorphic. &  \begin{equation*}
        \rho_{tr}:T_0R\to H^1(X,\Theta_{tr})
    \end{equation*}\tabularnewline
    \hline
    Family of $h$-deformations & $\Theta_h$ locally given by vector fields of the form :
    \begin{equation*}
        \sum\limits_k a_{\mathcal{L},k}(r,z^{(\alpha)},w^{(\alpha)})\frac{\partial}{\partial z^{(\alpha)}_k}
    \end{equation*}
    \begin{equation*}
        +\sum\limits_k b_{tr, k}(r,w^{(\alpha)})\frac{\partial}{\partial w^{(\alpha)}_k},
    \end{equation*}
    with $a_{\mathcal{L},k}$ holomorphic and $b_{tr, k}$ holomorphic. &  \begin{equation*}
        \rho_h:T_0R\to H^1(X,\Theta_h)
    \end{equation*} \tabularnewline 
    \hline
    Family of $f$-deformations & $T\F^{1,0}$ locally given by vector fields of the form :
    \begin{equation*}
        \sum\limits_k a_{\mathcal{L},k}(r,z^{(\alpha)},w^{(\alpha)})\frac{\partial}{\partial z^{(\alpha)}_k},
    \end{equation*} 
    with $a_{\mathcal{L},k}$ holomorphic. & \begin{equation*}
        \rho_f:T_0R\to H^1(X,T\F^{1,0})
    \end{equation*}\tabularnewline
    \hline
 \end{tabular}\\

Note that the multifoliate structure point of view gives rise to another description of infinitesimal deformations. We will not make the difference between all those descriptions, since it will be clear which description we adopt.

\subsection{Kuranishi's Theorem}

There exists for the deformations we defined, particular families called versal families of deformations. We have :

\begin{definition}
   The family $\mathcal{X}/R=(\mathcal{X},R,(X,\mathcal{F}))$ is said to be versal (resp. universal) if for any other family of deformation $\mathcal{X}'/R'=(\mathcal{X}',R',(X,\mathcal{F}))$, there exist a map between analytic spaces $\Theta: (R',0)\to (R,0)$ such that $\mathcal{X}'/R'$ and $\mathcal{X}_\Theta/R'$ are isomorphic as families of deformations, and furthermore, the differential of $\Theta$ at $0$ is unique (resp. $\Theta$ is unique). 
\end{definition}

When there exists a versal family of deformations, the corresponding Kodaira-Spencer map is an isomorphism (see \cite{NIC} 2.10).
We can also show (\cite{Completness}) that up to an isomorphism, the analytic space parametrizing such a family is unique. 

\begin{definition}
    An analytic space $R$ parametrizing a versal deformation family is called versal deformation space or simply Kuranishi space.
\end{definition}

The name Kuranishi space associated with such spaces comes from the profound theorem due to Kuranishi in the context of deformations of compact complex manifolds stating that for any compact complex manifold there exists a versal family of deformations of its complex structure. We have in fact :

\begin{theoreme}
    (Kuranishi)
    For any compact complex manifold $X_0$, there exists a versal family of deformations parametrized by an analytic space germ $(R,0)$. 

    More precisely, there is an open neighborhood $U$ of $0$ in $H^1(X_0,TX_0)$ and an analytic map $\Psi:U\to H^2(X_0,TX_0)$ such that $(R,0)$ is the germ at $0$ of $\Psi^{-1}(0)$.
\end{theoreme}

The map $\Psi$ is constructed from the Maurer-Cartan equation that guarantees the integrability of a given almost complex structure obtained as a small deformation of the almost complex structure of $X_0$.\\

For holomorphic foliations, Girbau, Haefliger, Nicolau and Sundararaman showed \cite{GiNi,GHS} there are Kuranishi spaces for all three types of deformation, obtained by means of a Kuranishi-type theorem, considering the sheaf of germs of infinitesimal transformations adapted to the type of deformation, and its spaces of cohomologies $H^1$ and $H^2$.

We denote $K^h$, $K^f$ and $K^{tr}$ Kuranishi spaces corresponding respectively to the deformations of the types $h$, $f$ and $tr$, and $\mathcal{X}^h/K^h$, $\mathcal{X}^f/K^f$ and $\mathcal{X}^{tr}/K^{tr}$ the corresponding versal families. 
Denoting $H^i$ any of the cohomology spaces concerned, we have the following immediate corollary due to Kodaira and Spencer :

\begin{corollaire}(Kodaira-Spencer)
    If $H^2=0$, the Kuranishi space is smooth and isomorphic to an open neighborhood of $0$ in $H^1$.
\end{corollaire}

In the case of Calabi-Yau manifolds, results of Tian, Todorov and Bogomolov allows to show that such a result remains true in the absence of the cohomological assumption of the previous corollary. 

One aim of this text is to obtain such a result for Calabi-Yau type foliations for the space $K^f$ and obtain a decomposition of $K^h$ as the product of $K^f$ and $K^{tr}$.

\section{Unobstruction of $f$-deformations for Strongly Calabi-Yau foliations}\label{unobstruction}

The purpose of this section is to show the formal existence of a solution to the Maurer-Cartan's equation for any infinitesimal $f$-deformation of any strongly Calabi-Yau foliation. The goal here is to show :

\begin{theoreme}\label{lissite}
Let $(X,\mathcal{F})$ be a compact Kähler manifold foliated by a strongly Calabi-Yau foliation and $v\in H^1(X,T\mathcal{F})$. Then, there exist a formal power serie $\sum\limits_{k>0}\phi_k t^k$, with $\phi_k\in \mathcal{A}^{0,1}(T\mathcal{F})$ such that :
\begin{equation*}
    \overline{\partial}\phi_1=0 \mbox{ and   } \ \overline{\partial}\phi_k=-\sum\limits_{i+j=k}[\phi_i,\phi_j]
\end{equation*}
with $[\phi_1]=v$ and such that for all $k$, $\Delta_\mathcal{F}(\phi_k)=0$.
\end{theoreme}

Firstly we define what strongly Calabi-Yau foliations are.

\subsection{Strongly Calabi-Yau foliations}

We can consider Calabi-Yau manifolds as compact Kähler manifolds on which there is a holomorphic global form of maximal degree which trivializes its canonical bundle $K_X$. If we adapt this for the foliated situation, we can at the beginning suppose that foliated manifold $(X,\mathcal{F})$ of Calabi-Yau type are compact Kähler manifolds on which there is a global form with a degree equal to the dimension of the foliation which trivializes $K_\mathcal{F}$, the canonical bundle of the foliation. Let make this more precise.\\

Let $(X,\mathcal{F})$ be a compact Kähler foliated manifold of dimension $n$, with $\mathcal{F}$ a regular holomorphic foliation of dimension $p$. We fix an atlas $(U_\alpha,(x_1^{(\alpha)},\dots,x_n^{(\alpha)}))_\alpha$ equal to $(U_\alpha,(z^{(\alpha)}_1,\dots,z^{(\alpha)}_p,w^{(\alpha)}_1,\dots,w^{(\alpha)}_q))_\alpha$, where $q=n-p$, adapted to the foliation. Since we are dealing with holomorphic foliations, we will simply denote $\TF$ rather than $T\F^{1,0}$ the holomorphic tangent bundle of $\mathcal{F}$ and $\nF$ its normal bundle. We have :

\begin{equation*}
    0 \xrightarrow[]{} \TF \xrightarrow[]{} TX \xrightarrow[]{} \nF \xrightarrow[]{} 0
\end{equation*}

and :

\begin{equation*}
    0 \xrightarrow[]{} \nF^* \xrightarrow[]{} TX^* \xrightarrow[]{} \TF^* \xrightarrow[]{} 0.
\end{equation*}

and then :

\begin{equation*}
    0 \xrightarrow[]{} N_k^* \xrightarrow[]{} \bigwedge^k TX^* \xrightarrow[]{} \bigwedge^k T\mathcal{F}^* \xrightarrow[]{} 0
\end{equation*}

Locally, $N_*^*$ is the ideal of $\Lambda^* TX^*$ generated by $<dw_1,\dots,dw_q>$. We finally have : 

\begin{equation*}
    0 \xrightarrow[]{} \bigwedge^l \overline{TX^*}\otimes N_k^* \xrightarrow[]{} \bigwedge^l \overline{TX^*}\otimes\bigwedge^k TX^*\xrightarrow[]{} \bigwedge^l \overline{TX^*}\otimes\bigwedge^k T\mathcal{F}^* \xrightarrow[]{} 0.
\end{equation*}

We denote by $\AqNps$, $\Apq$ and $\AqTFps$ the sheaves of germs of sections of each bundle of the precedent exact sequence and we define projections $\Pi$ for each pair $(k,l)$ (without noting them) by :

\begin{equation}\label{suiteexacteformes}
    0 \xrightarrow[]{} \AqNps \xrightarrow[]{} \Apq \xrightarrow[]{\Pi} \AqLpTFs \xrightarrow[]{} 0.
\end{equation}

We also denote $\AAqNps$, $\AApq$ et $\AAqTFps$ the spaces of global sections of each sheaf. We remark by immediate calculations :

\begin{prop}\label{propdrond}
    The differential operator $\partial$ preserves $\AqNss$ and defines 
    \begin{equation*}
        \pF:\AqTFps\to \AqTFpsplus 
    \end{equation*}
    such that $\pF\circ\Pi=\Pi\circ\partial$.
\end{prop}

We finally denote by $\AAqLpTX$ the global sections of $\AqLpTX$ and we define on $\bigoplus\limits_{k,l}A^{0,l}(\bigwedge^k TX)$ a supercommutative algebra structure via the $\mathbb{Z}/2\mathbb{Z}$-graduation given by $k+l$ modulo $2$. We have for $\alpha\in A^{0,l}_X$ and $v\in\Gamma(X,\bigwedge^k TX)$ a global $k$-vector fields, $\alpha\otimes v=(-1)^{kl}v\otimes \alpha$. 

For a Calabi-Yau manifold $Y$, the existence of a global holomorphic form $\Omega$ vanishing nowhere of maximal degree (and then trivializing $K_Y$) allows us to construct the following isomorphism : 
\begin{equation*}
\begin{array}{ccccc}
    \eta_Y & : & A^{0,l}_Y(\bigwedge^k TY) 
  & \simeq & A^{n-k,l}_Y\\
     & & \alpha_I\otimes v_J & \mapsto & (v_J \lrcorner \Omega) \wedge \alpha_I.
\end{array}
\end{equation*}

In foliated situation, it is the existence of a global holomorphic form $\Omega_\mathcal{F}\in A^{p,0}_X$, with $p=dim(\mathcal{F})$, such that $\Pi(\Omega_\mathcal{F})$ trivializes $K_\F:=\bigwedge^p T\mathcal{F}^*$ that allows the existence of the following isomorphism :

\begin{equation*}
\begin{array}{ccccc}
    \eta_{\mathcal{F}} & : & \AAqLpTF 
  & \simeq & \AAqLkmpTFs\\
     & & \alpha_I\otimes v_J & \mapsto & \Pi((v_J \lrcorner \Omega_\mathcal{F}) \wedge \alpha_I),
\end{array}
\end{equation*}
where $(v_1\wedge\dots\wedge v_k)\lrcorner\Omega_\mathcal{F}=i_{v_1}\circ\dots\circ i_{v_k}(\Omega_\mathcal{F})$, and $i_{v_i}$ contraction by $v_i$. We will write $(\alpha_I\otimes v_J)\lrcorner \Omega_{\mathcal{F}}$ rather than $(v_J \lrcorner \Omega_\mathcal{F}) \wedge \alpha_I$, understanding that contraction of $\Omega_\mathcal{F}$ by a holomorphic vector field gives a holomorphic form and that we respect classical notation for forms of any bidegree.

This hypothesis gives in particular a section of the precedent exact sequence as :

\begin{equation*}
    \begin{array}{ccccc}
         S & : & \AAqLpTFs & \to & A^{p-k,l}_X  \\
         & & \beta & \mapsto & \eta_\mathcal{F}^{-1}(\beta)\lrcorner\Omega_\mathcal{F}
    \end{array}
\end{equation*}

We give the following definition :

\begin{definition}
We call strongly Calabi-Yau foliation a regular holomorphic foliation $\mathcal{F}$ of dimension $p$ on a compact Kähler manifold $X$ endowed with a global holomorphic $p$-form $\Omega_\mathcal{F}\in H^0(X,\Omega_X^p)$, where $\Omega_X^p$ denotes the sheaf of holomorphic $p$-forms on $X$, such that $\Pi(\Omega_\mathcal{F})$ trivializes $K_\F$.
\end{definition}

For this type of foliations, we can construct the following commutative diagramm :

\begin{equation*}
   \begin{tikzcd}
 \AAqLpTF\arrow[r, "\eta_\mathcal{F}"] \arrow[rr, bend left , "\cdot\lrcorner\Omega_\mathcal{F} "] \arrow[d,"\Delta_\mathcal{F}"]
  & \AAqLkmpTFs\arrow[r, right =1ex,"S" ]\arrow[d,"\partial_\mathcal{F}"]
  & A^{p-k,l}_X\arrow[l, shift left=1ex, "\Pi"]\arrow[d,"\partial"]
  \\
A^{0,l}(\bigwedge^{k-1}\mathcal{TF})
& A^{0,l}(\bigwedge^{p-k+1}T\mathcal{F}^*)\arrow[l,"\eta^{-1}_\mathcal{F}"]
& A^{p-k+1,l}_X\arrow[l, "\Pi"] \\ 
\end{tikzcd} 
\end{equation*}

In particular we define $\Delta_\mathcal{F}$ as :

\begin{equation*}
    \begin{array}{ccccc}
         \Delta_\mathcal{F} & : & A^{0,l}(\bigwedge^{k}\mathcal{TF}) & \to & A^{0,l}(\bigwedge^{k-1}\mathcal{TF})  \\
         & & \phi & \mapsto & \eta_\mathcal{F}^{-1}\circ \partial_\mathcal{F} \circ \eta_\mathcal{F}(\phi).
    \end{array}
\end{equation*}

One of the aims of this text is to adapt the proof of the non-obstruction theorem of Calabi-Yau manifolds (particularly known as Bogomolov-Tian-Todorov theorem) given in \cite{HUY}, to the case of strongly Calabi-Yau foliations.

\begin{remarque}
    We can call Calabi-Yau foliation a holomorphic regular foliation on a compact Kähler manifold $X$ whose canonical bundle is trivial ($K_\F\simeq \mathcal{O}_X$). Actually, Loray, Pereira and Touzet showed in \cite{SING} that the notions of Calabi-Yau foliations and strongly Calabi-Yau foliations are the same in numerous situations (for example if $\F$ is Calabi-Yau it suffices that $\F$ has a compact leaf, or $K_X$ is psef). It is even conjectured that they are the same notions by Sommese at the end of \cite{SOM}. In particular, when a Calabi-Yau foliation $\F$ is such that it admits an everywhere transverse and complementary foliation, then $\F$ is strongly Calabi-Yau.
\end{remarque}

\subsection{Tian-Todorov lemma}

The fundamental result for the proof of unobstruction theorem for Calabi-Yau manifolds is the Tian-Todorov lemma which can be understood as a way to express the bracket of two "$\partial$-closed" vector forms as a "$\partial$-exact" vector form. 

Therefore, we naturally try to obtain an adapted version of this result. Moreover, since $\TF$ is involutive, it leads us to think that the projection appearing in $\Delta_\mathcal{F}$ deletes only terms of no interest for the following equality :

\begin{lemme}(Tian-Todorov)\label{Tian-Todorov}
    Let $\phi\in A^{0,Q}(\mathcal{TF})$ and $\psi\in A^{0,l}(\mathcal{TF})$. We have :
    \begin{equation*}
        (-1)^l[\phi,\psi]=\Delta_\mathcal{F}(\phi\wedge\psi)-(-1)^l\Delta_\mathcal{F}(\phi)\wedge \psi -\phi\wedge \Delta_\mathcal{F}(\psi)
    \end{equation*}
\end{lemme}

\begin{preuve}
    We prove this locally. Let $(z_1,\dots,z_p,w_1,\dots,w_{n-p})$ local coordinates adapted to $\mathcal{F}$. We will denote $x_i:=z_i$ and $x_{p+i}:=w_i$. We can then suppose without loss of generality that $\phi$ and $\psi$ are locally of the form $\phi=\phi_{I,i}d\overline{x_I}\otimes \frac{\partial}{\partial z_i}$
    and $\psi=\psi_{J,j}d\overline{x_J}\otimes \frac{\partial}{\partial z_j}$.

    Locally, $\Omega_\mathcal{F}$ is of the form $\Omega_\mathcal{F}=f dz_1\wedge \dots \wedge dz_p+\Omega'$, with $\Pi(\Omega')=0$, i.e $\Omega'\in <dw_1,\dots,dw_{n-p}>$, and $f$ holomorphic without zero.

    Then :
    \begin{equation*}
    \begin{split}
    \begin{array}{ccc}
       \phi\lrcorner\Omega_\mathcal{F} & = & (-1)^{i-1} \phi_{I,i}fdz_1\wedge\dots\wedge \hat{dz_i}\wedge \dots\wedge dz_p\wedge d\overline{x_I}
        +\phi\lrcorner\Omega'  \\
        \Pi \partial (\phi\lrcorner\Omega_\mathcal{F}) & = &  \frac{\partial(\phi_{I,i}f)}{\partial z_i}dz_1\wedge\dots\wedge  dz_p\wedge d\overline{x_I}+0
        \\
        \Delta_\mathcal{F}(\phi) & = & \eta_\mathcal{F}^{-1}(\Pi \partial (\phi\lrcorner\Omega_\mathcal{F})\\
        & = & \frac{1}{f}\frac{\partial(\phi_{I,i}f)}{\partial z_i} d\overline{x_I}
    \end{array}     
    \end{split}
    \end{equation*}
    and, 
    \begin{equation*}
    \begin{split}
        \begin{array}{ccc}
        \Delta_\mathcal{F}(\phi)\wedge\psi & = & \frac{1}{f}\frac{\partial(\phi_{I,i}f)}{\partial z_i}\psi_{J,j} (d\overline{x_I}\wedge d\overline{x_J})\otimes \frac{\partial}{\partial z_j}\\
    \end{array}
    \end{split}
    \end{equation*}
    Likewise,
    \begin{equation*}
    \begin{split}
    \begin{array}{ccc}
        \phi\wedge\Delta_\mathcal{F}(\psi) & = & \frac{1}{f}\phi_{I,i}\frac{\partial(\psi_{J,j}f)}{\partial z_j}d\overline{x_I}\otimes\frac{\partial}{\partial z_i}\wedge d\overline{x_J}\\
        & = & (-1)^l\frac{1}{f}\phi_{I,i}\frac{\partial(\psi_{J,j}f)}{\partial z_j}d\overline{x_I}\wedge d\overline{x_J}\otimes\frac{\partial}{\partial z_i}\\
    \end{array}
    \end{split}
    \end{equation*}
    Moreover, 
    \begin{equation*}
        \phi\wedge\psi=(-1)^l \phi_{I,i}\psi_{J,j}(d\overline{x_I}\wedge d\overline{x_J})\otimes (\frac{\partial}{\partial z_i}\wedge\frac{\partial}{\partial z_j})
    \end{equation*}
    and, without loss of generality, supposing $i<j$, 
    \begin{equation*}
    \begin{array}{ccc}
       (\phi\wedge\psi)\lrcorner\Omega_\mathcal{F} & = & (-1)^l(-1)^{i-1+j-1}\phi_{I,i}\psi_{J,j}fdz_1\wedge\dots\wedge \hat{dz_i}\wedge \dots\wedge\hat{dz_j}\wedge \dots\wedge dz_p\wedge d\overline{x_I}\wedge d\overline{x_J}\\
       & & + (-1)^l(\phi\wedge\psi)\lrcorner\Omega'\\
       \Pi\partial((\phi\wedge\psi)\lrcorner\Omega_\mathcal{F}) & = & (-1)^l(-1)^{j-1}\frac{\partial (\phi_{I,i}\psi_{J,j}f)}{\partial z_i}dz_1\wedge\dots\wedge\hat{dz_j}\wedge \dots\wedge dz_p\wedge d\overline{x_I}\wedge d\overline{x_J}\\
       & & + (-1)^l(-1)^{i}\frac{\partial (\phi_{I,i}\psi_{J,j}f)}{\partial z_j}dz_1\wedge\dots\wedge\hat{dz_i}\wedge \dots\wedge dz_p\wedge d\overline{x_I}\wedge d\overline{x_J}\\
       \Delta_{\mathcal{F}}(\phi\wedge\psi) & = & \frac{1}{f}(-1)^l\frac{\partial (\phi_{I,i}\psi_{J,j}f)}{\partial z_i} (d\overline{x_I}\wedge d\overline{x_J})\otimes\frac{\partial}{\partial z_j}\\
       & & -\frac{1}{f}(-1)^l\frac{\partial (\phi_{I,i}\psi_{J,j}f)}{\partial z_j} (d\overline{x_I}\wedge d\overline{x_J})\otimes\frac{\partial}{\partial z_i}\\
       
       & = &(-1)^l(\frac{1}{f}\psi_{J,j}\frac{\partial (\phi_{I,i}f)}{\partial z_i}+\phi_{I,i}\frac{\partial \psi_{J,j}}{\partial z_i}) (d\overline{x_I}\wedge d\overline{x_J})\otimes\frac{\partial}{\partial z_j}\\
       & & -(-1)^l(\frac{1}{f}\phi_{I,i}\frac{\partial (\psi_{J,j}f)}{\partial z_j}+\psi_{J,j}\frac{\partial \phi_{I,i}}{\partial z_j}) (d\overline{x_I}\wedge d\overline{x_J})\otimes\frac{\partial}{\partial z_i}\\
    \end{array}
    \end{equation*}
    Finally, we get :
    \begin{equation*}
        \begin{array}{ccc}
             \Delta_\mathcal{F}(\phi\wedge\psi)-(-1)^l\Delta_\mathcal{F}(\phi)\wedge\psi-\phi\wedge\Delta_\mathcal{F}(\psi) & = & (-1)^l(\phi_{I,i}\frac{\partial \psi_{J,j}}{\partial z_i} (d\overline{x_I}\wedge d\overline{x_J})\otimes\frac{\partial}{\partial z_j}\\
       & & -\psi_{J,j}\frac{\partial \phi_{I,i}}{\partial z_j} (d\overline{x_I}\wedge d\overline{x_J})\otimes\frac{\partial}{\partial z_i})\\
       & = & (-1)^q[\phi,\psi]
        \end{array}
    \end{equation*}

\end{preuve}

We immediatly get the following :

\begin{corollaire}
Let $\phi\in A^{0,Q}(\mathcal{TF})$, $\psi\in A^{0,l}(\mathcal{TF})$ such that $\Delta_\mathcal{F}(\phi)=0$ and $\Delta_\mathcal{F}(\psi)=0$.

Then, we get : 
\begin{equation*}    
(-1)^q[\phi,\psi]=\Delta_\mathcal{F}(\phi\wedge\psi) \mbox{ et } (-1)^q\eta_\mathcal{F}[\phi,\psi]=\partial_\mathcal{F} (\eta_\mathcal{F}(\phi\wedge\psi).
\end{equation*}
\end{corollaire}

\subsection{Commutation relations}

We show here that $\overline{\partial}$ commute with $\Pi$, and up to a sign $\eta_\mathcal{F}$. After the proposition we will denote by $\overline{\partial}$ any of the following : $\overline{\partial}$, $\overline{\partial}_{\bigwedge^* T\mathcal{F}}$ and $\overline{\partial}_{\bigwedge^* T\mathcal{F}^*}$.

\begin{prop}\label{propcommutation}
We have the following relations :
\begin{equation*}
\begin{array}{cccc}
     &\overline{\partial}_{\bigwedge^* T\mathcal{F}^*}\circ\Pi & = & \Pi\circ\overline{\partial}  \\
     \forall \beta\in A^{r,s}_X, \ \overline{\partial}\mbox{-closed }, &\overline{\partial}((\ \cdot \ ) \lrcorner \beta) & = & (-1)^{r-k}(\overline{\partial}_{\bigwedge^* T\mathcal{F}} ( \ \cdot \ ))\lrcorner\beta \\
     &\overline{\partial}_{\bigwedge^* T\mathcal{F}^*}\circ\eta_\mathcal{F} & = & (-1)^{p-k} \eta_\mathcal{F}\circ \overline{\partial}_{\bigwedge^* T\mathcal{F}}
\end{array}
\end{equation*}
\end{prop}

\begin{preuve}
    \underline{(i) :} Let $\alpha\in A^{k,m}_X$. Locally we have $\alpha=\sum\limits_{I\subset\{1,\dots,p\}^k} \alpha_{I,J}dx_I\wedge d\overline{x_J}+\sum\limits_{I\not\subset\{1,\dots,p\}^k} \alpha_{I,J}dx_I\wedge d\overline{x_J}$. In particular, $\Pi(\alpha)=\sum\limits_{I\subset\{1,\dots,p\}^k} \alpha_{I,J}dx_I\wedge d\overline{x_J}$, then we have the following equality, $\overline{\partial}_{\bigwedge^* T\mathcal{F}^*}\circ\Pi(\alpha)=\sum\limits_{I\subset\{1,\dots,p\}^k,l} \frac{\partial \alpha_{I,J}}{\partial \overline{x_l}}d\overline{x_l}\wedge dx_I\wedge d\overline{x_J}$. 

    Likewise, $\overline{\partial} \alpha=\sum\limits_{I\subset\{1,\dots,p\}^k,l}\frac{\partial \alpha_{I,J}}{\partial \overline{x_l}} d\overline{x_l}\wedge dx_I\wedge d\overline{x_J}+\sum\limits_{I\not\subset\{1,\dots,p\}^k,l} \frac{\partial \alpha_{I,J}}{\partial \overline{x_l}}d\overline{x_l}\wedge dx_I\wedge d\overline{x_J}$, and then, $\Pi\circ\overline{\partial}(\alpha)=\sum\limits_{I\subset\{1,\dots,k\}^p,l}\frac{\partial \alpha_{I,J}}{\partial \overline{x_l}} d\overline{x_l}\wedge dx_I\wedge d\overline{x_J}=\overline{\partial}_{\bigwedge^* T\mathcal{F}^*}\circ\Pi(\alpha)$.

    \underline{(ii) :} Let $\phi\in A^{0,m}(\bigwedge^{k}\mathcal{TF})$, which we suppose locally of the form $\phi=\alpha\otimes \frac{\partial}{\partial x_I}$. We have for any global $\overline{\partial}$-closed form $\beta=\sum\limits_{K,L}\beta_{K,L}dx_K\wedge d\overline{x_L}$ of bidegree $(r,s)$, the equality $\phi\lrcorner\beta=\sum\limits_{K \supset I,L}(-1)^{\sum\limits_{i\in I}(i-1)}\beta_{K,L}dx_{K\textbackslash I}\wedge\alpha\wedge d\overline{x_L}$. However, $\overline{\partial}\beta =0$ implies $\overline{\partial} (\sum\limits_{K\supset I, L}\beta_{K,L}d\overline{x_L})=0$. We get then $\overline{\partial}(\phi\lrcorner\beta)=\sum\limits_{K \supset I,L}(-1)^{\sum\limits_{i\in I}(i-1)}(-1)^{r-k}dx_{K\textbackslash I}\wedge\overline{\partial}(\alpha\wedge\beta_{K,L} d\overline{x_L})$, i.e, $\overline{\partial}(\phi\lrcorner\beta)=(-1)^{r-k}\sum\limits_{K \supset I,L}(-1)^{\sum\limits_{i\in I}(i-1)}dx_{K\textbackslash I}\wedge\overline{\partial}(\alpha)\wedge\beta_{K,L} d\overline{x_L}$, and then, $\overline{\partial}(\phi\lrcorner\beta)=(-1)^{r-k}\overline{\partial}(\phi)\lrcorner\beta$.

    \underline{(iii) :} We apply $(ii)$ to $\Omega_\mathcal{F}=fdx_1\wedge\dots\wedge dx_p+\Omega'$. The composition by $\Pi$ of $(ii)$ allows us to delete the term $\cdot\lrcorner \Omega'$ and we get what we want.
\end{preuve}

In the following, we will use those relations without mentioning them.

\subsection{Theorem proof}

We finish this section with the proof of theorem \ref{lissite}, which means that we can integrate any infinitesimal deformation of order $1$ for a strongly Calabi-Yau foliation.

\begin{preuve}

    \underline{Step 1 :} At first, we want to take a good representative of $v$. Let $\alpha_1\in \mathcal{A}^{0,1}(T\mathcal{F})$ be a representative of $v$. 

     Thanks to proposition $\ref{propcommutation}$ and the fact that $\alpha_1$ and $\Omega_\F$ are $\overline{\partial}$-closed, we get that $\alpha_1\lrcorner\Omega_\mathcal{F}\in\mathcal{A}^{p-1,1}(X)$ is still $\overline{\partial}$-closed. Then, we denote by $H$ its harmonic representative and we define $\phi_1:=\eta_\mathcal{F}^{-1}\circ\Pi(H)\in A^{0,1}(T\F)$.

    We have to verify that $\phi_1$ represents $v$ and that $\Delta_\mathcal{F}(\phi_1)=0$. 
     We know that $H=\alpha_1\lrcorner\Omega_\mathcal{F}+\overline{\partial}\beta$, with $\beta\in\mathcal{A}^{p-1,0}(X)$ thus thanks to commutation relations $\phi_1=\eta_\mathcal{F}^{-1}\circ\Pi(H)=\alpha_1+\overline{\partial}\beta '$, with $\beta'\in \mathcal{A}^{0,0}(T\mathcal{F})$, so $[\phi_1]=v$.

    Since $H$ is $\partial$-closed, we have $\partial_\mathcal{F}(\Pi(H))=0$ and then $\Delta_\mathcal{F}(\phi_1)=\eta_\mathcal{F}^{-1}\circ\partial_\mathcal{F}(\eta_\mathcal{F}(\phi_1))=0$.

    \underline{Step 2 :} We construct $\phi_2$ in this second step. 
    Thanks to the precedent step and the lemma \ref{Tian-Todorov} we get $-[\phi_1,\phi_1]=\Delta_\mathcal{F}(\phi_1\wedge\phi_1)=\eta_\mathcal{F}^{-1}\circ\Pi(\partial ((\phi_1\wedge\phi_1)\lrcorner\Omega_\mathcal{F}))$. We also have that $\partial ((\phi_1\wedge\phi_1)\lrcorner\Omega_\mathcal{F})$ is $\partial$-exact and is $\overline{\partial}$-closed thanks to proposition $\ref{propcommutation}$ since $\phi_1\wedge\phi_1$ is $\overline{\partial}$-closed. Thus, by $\partial\overline{\partial}$-lemma, we have $[\phi_1,\phi_1]=\eta_\mathcal{F}^{-1}\circ\Pi(\overline{\partial}\partial\alpha_2)$. Then we define $\phi_2:=(-1)^{p}\eta_\mathcal{F}^{-1}\circ\Pi(\partial\alpha_2)$ and verify $\overline{\partial}\phi_2=-\eta_{\mathcal{F}}^{-1}\circ\Pi(\overline{\partial}\partial \alpha_2)=-[\phi_1,\phi_1]$. Finally, by proposition $\ref{propdrond}$ applied to $\partial \alpha_2$, $\Delta_\mathcal{F}(\phi_2)=0$.
    
    \underline{Step 3 :} Suppose that $(\phi_i)_{1\leq i\leq k}$ have been constructed, we want to define $\phi_{k+1}$. The brackets $[\phi_i,\phi_j]$ are equal to $\Delta_\mathcal{F}(\phi_i\wedge\phi_j)=\eta_\mathcal{F}^{-1}\circ\Pi(\partial((\phi_i\wedge\phi_j)\lrcorner\Omega_\mathcal{F}))$, by foliated Tian-Todorov lemma and the hypothesis on the $\phi_i$. In the same time we have, 
    \begin{equation*}
    \begin{array}{ccc}
        \overline{\partial}\sum\limits_{1<i<k}(\phi_i\wedge\phi_{k-i}) & = & \sum\limits_{1<i<k}(\overline{\partial}\phi_i\wedge\phi_{k-i})-(\phi_i\wedge\overline{\partial}\phi_{k-i})  \\
         & = &  \sum\limits_{1<i<k}((\sum\limits_{1<j<i}[\phi_j,\phi_{i-j}]\wedge\phi_{k-i})-(\phi_i\wedge\sum\limits_{1<j<k-i}[\phi_j,\phi_{k-i-j}])) \\
         & = & \sum\limits_{1<i<k}(\sum\limits_{1<j<i}[\phi_j,\phi_{i-j}]\wedge\phi_{k-i})-\sum\limits_{1<i<k}(\phi_i\wedge\sum\limits_{1<j<k-i}[\phi_j,\phi_{k-i-j}]) \\
         & = & \sum\limits_{1<i<k}(\sum\limits_{1<j<i}[\phi_j,\phi_{i-j}]\wedge\phi_{k-i})-\sum\limits_{1<k-i<k}(\phi_{k-i}\wedge\sum\limits_{1<j<i}[\phi_j,\phi_{i-j}]) \\
         & = & \sum\limits_{1<i<k}(\sum\limits_{1<j<i}[\phi_j,\phi_{i-j}]\wedge\phi_{k-i}-\phi_{k-i}\wedge\sum\limits_{1<j<i}[\phi_j,\phi_{i-j}])\\
         & = & 0.
    \end{array}
    \end{equation*}

    Then proposition $\ref{propcommutation}$ implies that $\overline{\partial}\sum\limits_{1<i<k}((\phi_i\wedge\phi_{k-i})\lrcorner\Omega_\mathcal{F})=0$.
    
    Thus, $\partial(\sum\limits_{1<i<k}((\phi_i\wedge\phi_{k-i})\lrcorner\Omega_\mathcal{F}))$ is $\partial$-exact and $\overline{\partial}$-closed. We apply $\partial\overline{\partial}$-lemma and we get $\alpha_{k+1}$ such that $\partial(\sum\limits_{1<i<k}((\phi_i\wedge\phi_{k-i})\lrcorner\Omega_\mathcal{F}))=\overline{\partial}\partial \alpha_{k+1}$. Then we define $\phi_{k+1}:=(-1)^{k-1}\eta_{\mathcal{F}}^{-1}\circ\Pi(\partial \alpha_{k+1})$ and by the lemma we have $\Delta_\mathcal{F}(\phi_{k+1})=0$. 
\end{preuve}

The convergence of this formal serie in a neighborhood of $0$ is due to a theorem of Artin applied on $X$ without refering the foliation, furnishing an analytic solution to our problem of deformation.\\

We can rewrite our theorem as :

\begin{theorem}\label{lisse}
    Let $(X,\mathcal{F})$ be a compact Kähler manifold endowed with a strongly Calabi-Yau foliation. Then, the Kuranishi space $K^f$ of its $f$-deformation is smooth. 
\end{theorem}

\section{A theorem of decomposition}\label{dec}

In their article \cite{GiNi}, Girbau and Nicolau express the Kuranishi space $K^h$ as $K^f\times K^{tr}$, in the sense of analytic germ spaces, in the situation where the regular holomorphic foliation $(X,\mathcal{F})$ they deform admits a foliation everywhere transverse and complementary to $\F$, and with $H^2(X,T\F)~=~0$.

This last hypothesis allows them to construct a solution to the extension problem they consider in a very direct way.

In the case of a strongly Calabi-Yau foliation, we always have the existence of a foliation everywhere transverse and complementary to $\F$. To distinguish Kuranishi spaces of each foliation, from now on we will denote them $K_\F^h$, $K^f_\F$ and $K^{tr}_\F$.  The second hypothesis is no longer checked, but the vanishing of $H^2(X,T\mathcal{F})$ and the strongly Calabi-Yau hypothesis both imply that $K^f$ is smooth. 

We then want to check if the decomposition theorem is still true when the foliation is strongly Calabi-Yau, and better, if it is sufficient for $K^f_\F$ to be smooth (with the additional hypothesis about the existence of the transverse and complementary foliation). \\

So, the goal here is to show that the germs of space $K^f_\F\times K^{tr}_\F$ et $K^h_\F$ are isomorphic, under the hypothesis that $\F$ is a regular holomorphic foliation on $X$ admitting an everywhere transverse and complementary foliation $\mathcal{G}$ and such that $K_\F^f$ is smooth. A way to do it is to apply the following result from (\cite{GiNi}, 2.1) :

\begin{prop}\label{lemmeiso}
    Let $f:K\times K'\to K''$ and $g:K''\to K'$ morphisms between analytic spaces, with $K$ smooth, such that :
    \begin{enumerate}
        \item the tangent application $d_0 f$ is an isomorphism,
        \item the following diagram commutes :
        \begin{equation*}
   \begin{tikzcd}
  K\times K' \arrow[rr, "f"] \arrow[rd, "p_2"]
  & & K'' \arrow[ld,"g" ]\\
  & K' &  \\ 
\end{tikzcd} 
\end{equation*}
    \end{enumerate}
    Then, $f$ is an isomorphism of germs of analytic spaces.
\end{prop}

So what we want is to construct a map $K^f_\F\times K^{tr}_\F\to K^h_\F$, such that $(1)$ and $(2)$ are satisfied.\\

Girbau and Nicolau showed in (\cite{GiNi}, 3.3) that when there is a foliation transverse and supplementary to $\F$, we have an isomorphism of analytic spaces $K^f_\F\simeq K^0_\F:=\pi^{-1}(0)$, where $\pi:K^h_\F\to K^{tr}_\F$ is the corresponding versal map seen as the forgetful of the complex structures of the leaves. This means in particular that the versal familly $\mathcal{X}^f/K^f_\F$ of $f$-deformation can be represented by a trivial family of deformation from the transverse deformations view point. The representation of such a familly can then be made by local coordinates parameterized by $K^0_\F$ of the form $(U_\alpha,(z^{(\alpha)}_1(r),\dots,z^{(\alpha)}_p(r),w^{(\alpha)}_1(0),\dots, w^{(\alpha)}_q(0)))_{r\in K^0}$, adapted at $0$ to the foliations $\mathcal{F}$ and $\mathcal{G}$ everywhere transverse. 

Under the same hypothesis, they showed that the versal family $\mathcal{X}^{tr}/K^{tr}_\F$ of transverse deformation can be represented by a family of holomorphic deformation. In particular, we can represent this family by $(U_\alpha, (z^{(\alpha)}_1(0), \dots, z^{(\alpha)}_p(0), w^{(\alpha)}_1(r),\dots, w^{(\alpha)}_q(r))_{r\in K^{tr}})_\alpha.$

Following the Kodaira-Spencer procedure, we recalled earlier for \ref{representationfamille}, we can construct jet forms respecting the initial multifoliate structure corresponding to each of theses families. We denote $u^f(r_f)=(\phi^f(r_f),\eta^f(r_f))$ and $u^{tr}(r_{tr})=(\phi^{tr}(r_{tr}),\eta^{tr}(r_{tr}))$ the respective jet forms, where $r_f\in K^0_\F\simeq K^f_\F$ et $r_{tr}\in K^{tr}_\F$.

\begin{prop}
    The jet forms $u^f(r_f)$ and $u^{tr}(r_{tr})$ are such that $\phi^f(r_f)=\omega^f(r_f)+\overline{\omega^f(r_f)}$ and $\phi^{tr}(r_{tr})=\omega^{tr}(r_{tr})+\overline{\omega^{tr}(r_{tr})}$ are of the form :
    
        \begin{multline*}
             \omega^f(r_f)  =  \sum\limits b_{l,k}^{(\alpha)}(r_f) \frac{\partial}{\partial z^{(\alpha)}_l (r_f)}\otimes dz^{(\alpha)}_k (0)\\ + \sum\limits b_{l,p}^{(\alpha)}(r_f) \frac{\partial}{\partial z^{(\alpha)}_l (r)}\otimes dw^{(\alpha)}_p (0) +\sum\limits  \frac{\partial}{\partial w^{(\alpha)}_p (0)}\otimes dw^{(\alpha)}_p (0)
             \end{multline*}
             \begin{multline*}
             \omega^{tr}(r_{tr})  =  \sum\limits \frac{\partial}{\partial z^{(\alpha)}_k (0)}\otimes dz^{(\alpha)}_k (0) + \sum\limits b_{q,p}^{(\alpha)}(r_{tr}) \frac{\partial}{\partial w^{(\alpha)}_q (r_{tr})}\otimes dw^{(\alpha)}_p (0)
        \end{multline*}
    
\end{prop}

\begin{preuve}
    See the construction of the first term of the form jet from the datum of local coordinates recalled after theorem \ref{representationfamille}.
\end{preuve}

This construction allows us to construct a new family of deformation parametrized by $K^0_\F\times K^{tr}_\F$ described by a jet form $u^h(r_f, r_{tr})$ obtained from $\phi^h(r_f, r_{tr})$, a global vector $1$-form (by virtue of $TX=T\F\oplus T\mathcal{G}$), defined by $\phi^h(r_f, r_{tr})=\omega^h(r_f, r_{tr})+\overline{\omega^h(r_f, r_{tr})}$, with :

\begin{equation}\label{omegah}
  \begin{array}{cccc}
     \omega^h(r_f, r_{tr}) & = & & \sum\limits b_{l,k}^{(\alpha)}(r_f) \frac{\partial}{\partial z^{(\alpha)}_l (r_f)}\otimes dz^{(\alpha)}_k (0) + \sum\limits b_{l,p}^{(\alpha)}(r_f) \frac{\partial}{\partial z^{(\alpha)}_l (r_f)}\otimes dw^{(\alpha)}_p (0) \\
     & & + & \sum\limits b_{q,p}^{(\alpha)}(r_{tr}) \frac{\partial}{\partial w^{(\alpha)}_q (r_{tr})}\otimes dw^{(\alpha)}_p (0).
  \end{array}
\end{equation}

Reproducing the Kodaira-Spencer procedure used for \ref{representationfamille} it allows to construct a global vector $2$-form $\eta^h$ in such a way that $u^h:=(\phi^h,\eta^h)$ respects the multifoliate structure of $\F$ such that $[u^h,u^h]=0$. 

Then we obtain a holomorphic family of deformation of $\F$. Thanks to the versality of $K^h_\F$, there exist a map between analytic spaces $\overset{\sim}{\alpha}:K^0_\F\times K^{tr}_\F\to K^h_\F$ such that $(\mathcal{X}^f\times \mathcal{X}^{tr})/(K^f_\F\times K^{tr}_\F)$ and $\mathcal{X}^h_{\overset{\sim}{\alpha}}/(K^f_\F\times K^{tr}_\F)$ are isomorphic as families of deformation.\\

Let us now recall that Kodaira-Spencer map corresponding to each of these families $\rho_f$, $\rho_{tr}$ and $\rho_{h}$ ensure $\rho_f : T_0 K^0_\F\overset{\sim}{\to} H^1(X,T\mathcal{F})$, $\rho_{tr}:T_0 K^{tr}_\F \overset{\sim}{\to} H^1(X,\Theta_{tr})$ and $\rho_{h}:T_0 K^h_\F\overset{\sim}{\to} H^1(X,\Theta_h)$, and that the following diagramm commutes :

\begin{equation}
   \begin{tikzcd}
    T_0 K^0_\F \times T_0 K^{tr}_\F \arrow[r,"d_0 \overset{\sim}{\alpha}" ] \arrow[rd, "\rho_h' "]
  & T_0 K^h_\F \arrow[d, "\rho_h"]\\
  & H^1(X,\Theta_h)
\end{tikzcd} 
\end{equation}

where $\rho_h '$ is the Kodaira-Spencer map for the family of $h$-deformations represented by $u^h$ and parameterized by $K^0_\F\times K^{tr}_\F$.\\

We also recall that a representative of any of these Kodaira-Spencer classes can be obtained from the description of each family by jet forms in the following way : if $v=(\phi,\eta)$ represents a family of $f$-deformation (respectively $tr$,$h$-deformation) and $\frac{\partial}{\partial r}\in T_0 K^0_\F$ (respectively $T_0 K^{tr}_\F$, $T_0 K^h_\F$), the class $\rho^f (\frac{\partial}{\partial r})$ (respectively $\rho^{tr} (\frac{\partial}{\partial t})$, $\rho^h (\frac{\partial}{\partial t})$) is represented by $\frac{\partial v}{\partial r}:=(\frac{\partial \phi}{\partial r},\frac{\partial \eta}{\partial r})=(\frac{\partial \phi}{\partial r},d\frac{\partial \phi}{\partial r})$ (see \cite{MUL}, 8.3).\\

Note that by definition of $u^{tr}$, the Kodaira-Spencer map $\rho^{tr}$ can be seen as taking values in a subspace $T_r$ of $H^1(X,\Theta_h)$ such that $T_r\simeq H^1(X,\Theta_{tr})$. In other words, $T_r$ is the image of a map in cohomology induced by a section $s$ (given thanks to $\Theta_h\subset TX=T\F\oplus T\mathcal{G}$) :

\begin{equation*}
  \begin{tikzcd}
  0 \arrow[r] 
  & T\mathcal{F} \arrow[r,"" ]
  & \Theta_h \arrow[r]
  & \Theta_{tr} \arrow[l, bend right , "s " above]\arrow[r]
  & 0.
\end{tikzcd} 
\end{equation*}


Then we show :

\begin{prop}
    The following diagramm commutes :
    \begin{equation}
   \begin{tikzcd}
    T_0 K^0_\F \times T_0 K^{tr}_\F \arrow[d,"\rho_f\times \rho_{tr}" left ] \arrow[rd, "\rho_h' "]\\
    H^1(X,T\mathcal{F})\oplus T_r \arrow[r,"S"] & H^1(X,\Theta_h)
\end{tikzcd} 
\end{equation}
and $S$, the sum, is an isomorphism between vector spaces.
\end{prop}

\begin{preuve}
    Let $\frac{\partial}{\partial r_f}\in T_0 K^0_\F$ and $\frac{\partial}{\partial r_{tr}}\in T_0 K^{tr}_\F$. We get :
    \begin{equation*}
        \rho_h ' (\frac{\partial}{\partial r_f},\frac{\partial}{\partial r_{tr}})=(\frac{\partial}{\partial r_f}+\frac{\partial}{\partial r_{tr}})(u_h).
    \end{equation*}
    An immediate calculation shows thanks to (\ref{omegah}) :
    \begin{equation*}
        (\frac{\partial}{\partial r_f}+\frac{\partial}{\partial r_{tr}})(u^h)=\frac{\partial u ^f}{\partial r_f}+\frac{\partial u^{tr}}{\partial r_{tr}}
    \end{equation*}
    and then :
    \begin{equation*}
        \rho_h ' (\frac{\partial}{\partial r_f}+\frac{\partial}{\partial r_{tr}})=\rho_f(\frac{\partial}{\partial r_f})+\rho_{tr}(\frac{\partial}{\partial r_{tr}})=:S(\rho_f(\frac{\partial }{\partial r_f}),\rho_{tr}(\frac{\partial}{\partial r_{tr}}))
    \end{equation*}
So the diagramm commutes.

Finally, the fact that $S$ is an isomorphism comes from the following exact sequence :
\begin{equation*}
  \begin{tikzcd}
  0 \arrow[r] 
  & T\mathcal{F} \arrow[r,"" ]
  & \Theta_h \arrow[r]
  & \Theta_{tr} \arrow[l, bend right , "s " above]\arrow[r]
  & 0.
\end{tikzcd} 
\end{equation*}
which is split.
\end{preuve}

From (2) and (3), we get in particular :

\begin{prop}
    The following diagramm commutes :
    \begin{equation}
   \begin{tikzcd}
    T_0 K^0_\F \times T_0 K^{tr}_\F \arrow[r,"d_0\overset{\sim}{\alpha}"] \arrow[d,"\rho_f\times \rho_{tr}" left ]  & T_0 K^h_\F\arrow[d, "\rho_h"]\\
    H^1(X,T\mathcal{F})\oplus T_r \arrow[r,"S"] & H^1(X,\Theta_h)
\end{tikzcd} 
\end{equation}

and therefore, $d_0\overset{\sim}{\alpha}$ is an isomorphism.
\end{prop}

We finaly get :

\begin{theorem}
    Let $K^f_\F$, $K^{tr}_\F$ and $K^h_\F$, the Kuranishi spaces corresponding respectively to the $f,tr,h$-deformations of a regular holomorphic foliation $(X,\mathcal{F})$ admitting a regular holomorphic foliation supplementary to $\mathcal{F}$ and such that $K^f_\F$ is smooth.

    Then there exists an isomorphism between analytic spaces $\alpha:K^f_\F\times K^{tr}_\F\to K^h_\F$ making commutative the following diagramm :
    \begin{equation*}
   \begin{tikzcd}
  K^f_\F\times K^{tr}_\F \arrow[rr, "\alpha"] \arrow[rd, "p_2"]
  & & K^h_\F \arrow[ld,"\pi" ]\\
  & K^{tr}_\F &  \\ 
\end{tikzcd} 
\end{equation*}
\end{theorem}

The following proof is an adaptation of a detailed proof of a result from (\cite{GiNi}, 2.3) that M. Nicolau kindly communicated to me.

\begin{preuve}
    According to the previous property, we are tempted to check that the diagramm :
    \begin{equation*}
   \begin{tikzcd}
  K^0_\F\times K^{tr}_\F \arrow[rr, "\overset{\sim}{\alpha}"] \arrow[rd, "p_2"]
  & & K^h_\F \arrow[ld,"\pi" ]\\
  & K^{tr}_\F &  \\ 
\end{tikzcd} 
\end{equation*}
commutes, which would suffice to apply the proposition \ref{lemmeiso} allowing us to conclude. The map $\overset{\sim}{\alpha}$ however has no reason to verify the previous relationship. To correct this, we construct a new $\alpha$ from $\overset{\sim}{\alpha}$.

We denote $\iota_{K^{tr}}:K^{tr}_\F\to K^0_\F\times K^{tr}_\F$ the map defined by $\iota_{K^{tr}}(r_{tr})=(0,r_{tr})$. Thus, the map $\Gamma:=\pi\circ \overset{\sim}{\alpha}\circ \iota_{K^{tr}}:K^{tr}_\F\to K^{tr}_\F$ induces a family $\mathcal{X}^{tr}_\Gamma /K^{tr}$ isomorphic to $\mathcal{X}^{tr}/K^{tr}$. By versality, $d_0\Gamma = id$ and therefore $\Gamma$ is an isomorphism, according to \cite{Completness}.

In addition we note by definition that $\overset{\sim}{\alpha}$ sends isomorphically $K^0_\F\times \{0\}$ onto $\pi^{-1}(0)=K^0_\F$.

We then denote $\alpha ':= \overset{\sim}{\alpha}\circ (id_{K^0}\times \Gamma^{-1}):K^0_\F\times K^{tr}_\F\to K^h_\F$ and $\Gamma':=\pi\circ \alpha '\circ \iota_{K^{tr}}:K^{tr}_\F\to K^{tr}_\F$. The map $\alpha '$ is still such that $d_0\alpha '$ is an isomorphism and sends isomorphically $K^0_\F\times \{0\}$ onto $\pi^{-1}(0)=K^0_\F$. The map $\Gamma '$ itself is such that :
\begin{equation*}
    \begin{array}{ccc}
         \Gamma ' (r_{tr}) & = & \pi\circ \alpha' (0,r_{tr}) \\
         & = & \pi\circ \overset{\sim}{\alpha}(0,\Gamma^{-1}(r_{tr}))  \\
         & = & \pi\circ \overset{\sim}{\alpha}\circ \iota_{K^{tr}}(\Gamma^{-1}(r_{tr}))\\
         & = & r_{tr},
    \end{array}
\end{equation*}
and therefore, $\Gamma ' =id_{K^{tr}}.$

We then denote, $h:=\Pi\circ \alpha ': K^0_\F\times K^{tr}_\F\to K^{tr}_\F$ and we define :

\begin{equation*}
    \begin{array}{ccccc}
H & : & K^0_\F\times K^{tr}_\F & \to & K^0_\F\times K^{tr}_\F \\
 & & (r_f,r_{tr}) & \mapsto & (r_f,h(r_f,r_{tr})) \\
\end{array}.
\end{equation*}

We note that $H(0,r_{tr})=(0,h(0,r_{tr}))=(0,r_{tr})$ and $H(r_f,0)=(r_f,h(r_f,0))=(r_f,0)$.
In particular, $d_0 H=id_{K^0_\F\times K^{tr}_\F}$ and therefore, $H$ is an automorphism of $K^0_\F\times K^{tr}_\F$.

We therefore define $G:=h\circ H^{-1}:K^0\times K^{tr}\to K^{tr}$. Thus, $G\circ H=h$ and therefore :
\begin{equation*}
    h(r_f,r_{tr})=G\circ H(r_f,r_{tr})=G(r_f,h(r_f,r_{tr})).
\end{equation*}

And from the fact that $h$ is surjective, we get $G=p_2$. We define $\alpha:=\alpha ' \circ H^{-1}$ and therefore :

\begin{equation*}
    \pi\circ \alpha=\pi\circ \alpha ' \circ H^{-1}=h \circ H^{-1}=G=p_2.
\end{equation*}

In addition, $d_0\alpha$ is still an isomorphism. Since $K^0_\F$ is smooth and $K^0_\F\simeq K^f_\F$, we get the result according to the proposition $\ref{lemmeiso}$.

\end{preuve}

To conclude this part, we want to show the following result :

\begin{equation*}
       K^{tr}_\F\simeq K^{0}_\mathcal{G}\simeq K^f_\mathcal{G}.
   \end{equation*}

for $(X,\F,\mathcal{G})$ two supplementary regular holomorphic foliations (in particular for $\F$ strongly Calabi-Yau foliation).

This result follows from the fact that $T^{0,1}=T\F^{0,1}\oplus T\mathcal{G}^{0,1}$ and therefore a small $f$-deformation of $\mathcal{G}$ can be represented by a global $1$-vector form corresponding to $\psi : T\mathcal{G}^{0,1}\to T\mathcal{G}^{1,0}$ (see the table in \ref{deformationtheory}). To see this, consider $\mathcal{X}^f_\mathcal{G}/K^f_\mathcal{G}$ the versal family given locally by $(z^{\alpha}_1(0),\dots,z^{\alpha}_p(0),w^{\alpha}_1(r),\dots, w^{\alpha}_q(r))_{\alpha, K^{f}_\mathcal{G}}$, such that for $r=0$ the coordinates are adapted to $T\F\oplus T\mathcal{G}=TX$ (see \cite{GiNi}), and so $T\F^{0,1}$ remains the same during deformation.

The subbundles $T^{0,1}$, $F$ and $G$ corresponding to the data of $\F$ and $\mathcal{G}$ are thus transformed in $\mathcal{X}^f_\mathcal{G}/K^f_\mathcal{G}$ into $T^{0, 1}_r=T\F^{0,1}\oplus T\mathcal{G}^{0,1}_r$, $F_r=T\F^{1,0}\oplus T\F^{0,1}\oplus T\mathcal{G}^{0,1}_r$ and $G_r=G$ respectively. In particular, $F_r$ is involutive and $G$ is involutive, so $T^{0,1}_r=F_r\cap G$ is involutive, and $T^{0,1}_r\subset F_r$ and $T^{0,1}_r\subset G$, and so $(T^{0,1}_r,F_r,G)$ corresponds to the data of two supplementary holomorphic foliations. In particular, we obtain :

\begin{prop}\label{rep multfi}
    The versal family of $f$-deformations of $\mathcal{G}$ can be represented by a family of deformations of the multifoliate structure $(X,\F,\mathcal{G})$.
\end{prop}

Forgetting the $\mathcal{G}$ foliation and the complex structure of the leaves, $\mathcal{X}^f_\mathcal{G} /K^f_\mathcal{G}$ can be seen as a family of $tr$-deformations of the transversely holomorphic foliation given by $F$. By versality, there exists $\beta:K^f_\mathcal{G}\to K^{tr}_\F$, such that $d_0\beta:H^1(X,T\mathcal{G})\to H^1(X,\Theta^{tr}_\F)$ is unique. Moreover, $\beta$ is such that $F_\beta(r)=F_r$ (from the versality porperty we have $\mathcal{X}^f_\mathcal{G} /K^f_\mathcal{G} \simeq (\mathcal{X}^{tr}_\F)_\beta / K^{f}_\mathcal{G}$).\\

Reciprocally, we can in this context obtain a holomorphic versal family of $tr$-deformations of $\F$ which we denote $\mathcal{X}^{tr}_\F/K^{tr}_\F$ locally given by coordinates of the particular precedent form $(z^{\alpha}_1(0), \dots,z^{\alpha}_p(0),w^{\alpha}_1(r'),\dots, w^{\alpha}_q(r'))_{\alpha, K^{tr}_\F}$, such that for $r'=0$ the coordinates are adapted to $T\F\oplus T\mathcal{G}=TX$ (see \cite{GiNi} 3.3).

In particular, since the family is holomorphic, the subbundles $T^{0,1}$ and $F$ are transformed into $T^{0,1}_{r'}$ and $F_{r'}$, so that $T\F^{0,1}_{r'}=T\F^{0,1}$, $T^{0,1}_{r'}$ and $F_{r'}$ are involutive, and $T^{0,1}_{r'}\subset F_{r'}$. We therefore have $F_{r'}=T\F^{0,1}\oplus T\F^{1,0}\oplus T\mathcal{G}^{0,1}_{r'}$, and the only transformations we are interested in are those corresponding to $1$-form vectors (global) of the form $\psi : T\mathcal{G}^{0,1}\to T\mathcal{G}^{1,0}$, since any other transformation ($\psi:T\mathcal{G}^{0,1}\to T\F^{0,1}\oplus T\F^{1,0}$) is offset by the fact that $T\F^{0,1}\oplus T\F^{1,0}\subset F_{r'}$. 

It follows that $G_{r'}=T\mathcal{G}^{0,1}_{r'}\oplus T\mathcal{G}^{1,0}_{r'}\oplus T\F^{0,1}=G$, and in particular, $G_{r'}=G$ is involutive, and we have $T^{0,1}_{r'}\subset G_{r'}=G$. In other words, $(T^{0,1}_{r'},F_{r'},G)$ corresponds, as before, to a small deformation of the multifoliate structure $(T^{0,1},F,G)$. Forgetting the foliation $\F$, we also note that $(T^{0,1}_{r'},G)$ is a $f$-deformation of the foliation $\mathcal{G}$. By versality, there exists an application $\gamma : K^{tr}_\F\to K^f_\mathcal{G}$ such that $d_0\gamma:H^1(X,\Theta^{tr}_\mathcal{F})\to H^1(X,T\mathcal{G})$ is unique. Moreover, $\gamma$ is such that $F_{\gamma(r')}=F_{r'}$ (from $(\mathcal{X}^f_\mathcal{G})_\gamma /K^f_\mathcal{G} \simeq \mathcal{X}^{tr}_\F / K^{f}_\mathcal{G}$).\\

In particular, for $r\in K^f_\mathcal{G}$, $F_{\gamma\circ\beta (r)}=F_{\beta(r)}=F_r$ and thus $\mathcal{X}^f_\mathcal{G}/K^f_\mathcal{G}\simeq (\mathcal{X}^f_\mathcal{G})_{\gamma\circ\beta}/K^f_\mathcal{G}$. Thus, by versality of $K^f_\mathcal{G}$, $\gamma\circ\beta$ is an isomorphism.

We then define $B:=\beta\circ (\gamma\circ\beta)^{-1}:K^f_\mathcal{G}\to K^f_\mathcal{G}$ and it follows that $B\circ\gamma=id:K^f_\mathcal{G}\to K^f_\mathcal{G}$. To conclude, according to \cite{Obstruction}, it suffices to show that $H^1(X,\Theta^{tr}_\mathcal{F})\simeq H^1(X,T\mathcal{G})$. Indeed Wavrik showed :

\begin{lemme}(Wavrik)
    Let $(S,0)$ and $(T,0)$ be germs of analytic spaces having the same tangential dimension. Let $\beta_1:(S,0)\to (T,0)$ and $\beta_2:(T,0)\to (S,0)$ be morphisms such that $\beta_2\circ \beta_1=id$. Then $\beta_2$ is an isomorphism.
\end{lemme}

So we prove the following to conclude :

\begin{lemme}
    Let $(X,\mathcal{F},\mathcal{G})$ be a compact complex manifold foliated by two supplementary regular holomorphic foliations (everywhere transverse and such that $TX=T\F\oplus T\mathcal{G}$). We then have :
    \begin{equation*}
        H^1(X,\Theta^{tr}_\F)\simeq H^1(X,T\mathcal{G}).
    \end{equation*}
\end{lemme}

\begin{preuve}
    We recall that :
    \begin{equation*}
        0\xrightarrow[]{}T\F\xrightarrow[]{} \Theta^h_\F\xrightarrow[]{} \Theta^{tr}_\F\xrightarrow[]{} 0,
    \end{equation*}
is a split exact sequence.
    So, in particular :
    \begin{equation*}
        0\xrightarrow[]{}H^1(X,T\F)\xrightarrow[]{} H^1(X,\Theta^h_\F) \xrightarrow[]{} H^1(X,\Theta^{tr}_\F)\xrightarrow[]{} 0.
    \end{equation*}
    And so, $H^1(X,\Theta^{tr}_\F)\simeq H^1(X,\Theta^h_\F)/H^1(X,T\F)$. Furthermore, the kernel of the projection $H^1(X,\Theta^{h}_\F)\to H^1(X,T\mathcal{G})$ is equal to $H^1(X,T\F)$, and so $H^1(X,\Theta^{tr}_\F) \xhookrightarrow{} H^1(X,T\mathcal{G})$ is injective.

    It remains to show that there exists an injective map $H^1(X,T\mathcal{G})\to H^1(X,\Theta^{tr}_\F)$.

    Let us note from Proposition \ref{rep multfi} that there exists a versal family of $f$-deformations of $\mathcal{G}$ locally given by coordinates $(U_\alpha, (z^{\alpha}_1(0),\dots,z^\alpha_p(0),w^\alpha_1(r),\dots,w^\alpha_q(r)))_{\alpha,K^f_\mathcal{G}}$ such that :
    \begin{equation*}
        \forall r, \alpha, \beta, i, j, \ \ \frac{\partial w^\alpha_i(r)}{\partial z^\beta_j(0)}=0.
    \end{equation*}
    It follows that the infinitesimal deformations corresponding to this family are locally of the form $\sum f_k(w)\frac{\partial}{\partial w^\alpha_k(0)}$. We therefore have an application $H^1(X,T\mathcal{G})\to H^1(X,\Theta^{tr}_\F)$ injective, which concludes the proof.
\end{preuve}

In conclusion, we have :

\begin{theorem}\label{KtrKf}
   Let $(X,\F,\mathcal{G})$ be two supplementary regular holomorphic foliations on $X$ a compact complex manifold. We have :
   \begin{equation*}
       K^{tr}_\F\simeq K^{0}_\mathcal{G}\simeq K^f_\mathcal{G}.
   \end{equation*}
   In particular, for $(X,\F)$ a strongly Calabi-Yau foliation and $\mathcal{G}$ a foliation supplementary to it, we have :
   \begin{equation*}
       K^{tr}_\F\simeq K^{0}_\mathcal{G}\simeq K^f_\mathcal{G}.
   \end{equation*}
\end{theorem}

If one ask me how should be named this result, I would say the Theorem of the sardines tin, or more briefly, Sardines theorem. In idea, a $tr$-deformation of $\F$ is about the deformation of a unique leaf of $\mathcal{G}$ while an $f$-deformation of $\mathcal{G}$ is about the deformation of all leaves of $\mathcal{G}$. One sardine moves, all sardine moves (be careful about the fact that the tin can be deformed too).

\section{Examples}

F. Touzet's classification \cite{CLATOU} provides examples of codimension 1 strongly Calabi-Yau foliations. We propose here to give the corresponding Kuranishi spaces for some examples in this classification.

\subsection*{i) $X$ is the product of a torus $T$ by a Calabi-Yau manifold $V$ (simply connected) and $\F$ is obtained by pulling back a foliation $\overset{\sim}{\F}$ linear of codimension 1 on the first factor by the canonical projection}

We have $X=T\times V\xrightarrow[]{p_1}T$ and $\F$ is given as $p_1^*\overset{\sim}{\F}$, with $dim(T)=n_T$, $dim(V)=n_V$, $n:=n_T+n_V$ et $\overset{\sim}{\mathcal{F}}$ of codimension $1$. In this situation, the existence of a foliation $\mathcal{G}$ everywhere transverse and complementary to $\F$ can be established by means of the wedge product of the pullback of the differential form defining a foliation transverse and complementary to $\overset{\sim}{\F}$ on $T$ by a holomorphic volume form of $V$.

Moreover, since $X$ is canonically trivial, then $K_\mathcal{G}$ is also trivial (by adjunction formula), so $\mathcal{G}$ is a Calabi-Yau foliation and is therefore strongly Calabi-Yau due to the existence of $\F$.

Finally, from Theorem \ref{KtrKf} we have $K^{tr}_\F \simeq K^f_\mathcal{G}$ , and so, according to the theorem \ref{lisse} for $f$-deformations of a strongly Calabi-Yau foliation, we have $K^{tr}_\F \simeq  U_\mathcal{G}$ , with $U_\mathcal{G} \subset H^1(X, T \mathcal{G})$ open.
Similarly, we have $K^f_\F \simeq U_\F$ , with $U_\F \subset H^1(X, T \F)$ open.

Finally, the decomposition theorem allows us to write :
\begin{equation*}
    K^h_\F\simeq K^{tr}_\F\times K^f_\F\simeq U_\mathcal{G}\times U_\F\subset H^1(X,T\mathcal{G})\oplus H^1(X,T\F), 
\end{equation*}

the last equality coming from the fact that $T X = T \mathcal{G} \oplus T \F$.

In particular, we note that $K^h_\F$ is an open of $H^1(X, T X)$, so it is also the Kuranishi space for the deformations of the complex structure of $X$.\\

We then try to make the link between the deformations of $\F$ and those of $\overset{\sim}{\F}$.\\

For $tr$-deformation, note that we have $\Theta_{\F}^{tr}=p_1^{-1}\Theta^{tr}_{\overset{\sim}{\F}}$ and $p_{1*} \Theta^{tr}_{\F}=\Theta^{tr}_{\overset{\sim}{\F}}$. Infinitesimaly we can say thanks to Leray spectral sequence and the exact sequence of low-degrees :

\begin{equation*}
    0\xrightarrow[]{} H^1(T,\Theta^{tr}_{\overset{\sim}{\F}})\xrightarrow[]{} H^1(X,\Theta^{tr}_\F)\to H^0(T,(\mathcal{R}^1p_{1*}\underline{\C})\otimes\Theta^{tr}_\F)=0,
\end{equation*}

since $V$ is simply connected.

In particular, we get $H^1(T,\Theta_{\overset{\sim}{\F}}^{tr})\simeq H^1(X,\Theta^{tr}_\F)$. Moreover, since each $tr$-deformation of $\overset{\sim}{\F}$ induces canonically a $tr$-deformation of $\F$, we can from the versal family $\mathcal{X}^{tr}_{\overset{\sim}{\F}}/K^{tr}_{\overset{\sim}{\F}}$ construct a family $\mathcal{X}/K^{tr}_{\overset{\sim}{\F}}$ of $tr$-deformations of $\F$ and obtain a morphism such that $K^{tr}_{\overset{\sim}{\F}}\hookrightarrow K^{tr}_\F$. Finally, thanks to the lemma \ref{lemmeiso}, we get : $K^{tr}_{\overset{\sim}{\F}}\simeq K^{tr}_\F$.\\

For $f$-deformation, note that $T\F=p_1^*T\overset{\sim}{\F}\oplus p_2^* TV$, and then $p_{1*}T\F=T\overset{\sim}{\F}\oplus (\mathcal{O}_T\otimes H^0(V,TV))$, so finaly $\mathcal{R}^1p_{1*}T\F=(\mathcal{R}^1p_{1*}p_1^*\mathcal{O}_X\otimes T\overset{\sim}{\F})\oplus (\mathcal{O}_T\otimes H^1(V,TV))$. But we have $\mathcal{R}^1p_{1*}p_1^*\mathcal{O}_X=0$ since $H^1(V,\mathcal{O}_V)=0$. In conclusion we have $\mathcal{R}^1p_{1*}T\F=\mathcal{O}_T\otimes H^0(V,TV)$. The exact sequence of low-degrees gives then :

\begin{multline*}
    0\xrightarrow[]{} H^1(T, T\overset{\sim}{\F})\oplus (H^1(T,\mathcal{O}_T)\otimes H^0(V,TV))\xrightarrow[]{} H^1(X,T\F)\\ \xrightarrow[]{} H^0(T,\mathcal{O}_T)\otimes H^1(V,TV)\xrightarrow{} H^2(T,p_{1*} T\F) 
\end{multline*}

On the one hand, we recognize $H^1(T, T\overset{\sim}{\F})\oplus (H^1(T,\mathcal{O}_T)\otimes H^0(V,TV))$ as the infinitesimal deformations generating deformations that preserves the product structure, i.e, we deform $(T,\overset{\sim}{\F})$ and $V$ independantly to obtain a family of deformations of $(T\times V,\F)$. On the other hand, we recognize $H^0(T,\mathcal{O}_T)\otimes H^1(V,TV)$ as the datum of a holomorphic family infinitesimal deformations of the complex structure of each $V$ for $p_1$.\\

Since any deformation of $(T,\overset{\sim}{\F})$ and of $V$ canonically induce deformations of $(X,\F)$, which are independant from each other, we can construct a morphism such that $K^f_{\overset{\sim}{\F}}\times K_V\hookrightarrow K^f_\F$ (which correspond to the deformations as a product).

\subsection*{ii) $X$ is a rational fibration over a manifold with a zero first Chern class
class (and hence described by Bogomolov's theorem) and the foliation $\F$ is transverse to the fibers}

The foliation $\mathcal{G}$ supplementary to $\F$ is given by the rational fibration
$p : X \to V$ (with $c_1(V ) = 0$ and for all $v \in V$ , $p^{-1}(v) = \mathbb{CP}^1$). Let $v_0 \in V$ be a point
of $V$.
From the decomposition theorem, it follows that :
\begin{equation*}
    K^h_\F\simeq K^f_\F\times K^{tr}_\F.
\end{equation*}

We partially recover the result of \cite{GiNi} (theorem 4.2):
\begin{theoreme}(Girbau \& Nicolau)
If the fibres $P$ of $p$ are Kähler and simply connected, then $K^f_\F\simeq K_V$, where $K_V$ is the Kuranishi space associated with the deformations of the complex structure of $V$, and furthermore, $K^h_\F \simeq K^f_\F \times K^{tr}_\F$.
\end{theoreme}

For the concrete data of $K^f_\F$ we apply this theorem, since $\mathbb{CP}^1$ is simply
connected.
Thus, $K^f_\F \simeq K_V \simeq U_V \subset H^1(V, T V )$, with $U_V$ open, according to Bogomolov-Tian-Todorov's theorem.

For the calculation of $K^{tr}_\F$ , Girbau, Haefliger and Sundaraman showed in (\cite{GHS}, 2.2) that $K^{tr}_\F$ is isomorphic to the versel space $K^{eq}_{\mathbb{CP}^1}$ of  $\pi_1(V, v_0)$-equivariant deformations of $\mathbb{CP}^1$. These deformations represent the way of transforming the representation $H : \pi_1 (V, v_0) \to Aut_\C(\mathbb{CP}^1)$ via the holonomy group associated with $(X, \F)$ (see
\cite{NIC}). 

More precisely, $K^{eq}_{\mathbb{CP}^1}$ is the analytic subspace of $(Aut_\C(\mathbb{CP}^1))^k$ ($Aut_\C(\mathbb{CP}^1$)
is a complex Lie group) given as the intersection of the analytic subspaces
$R_i(g_1,\dots , g_k) = 1$ and $\Sigma$, where $\pi_1 (V, v_0) =< g_1,\dots, g_k | R_1,\dots, R_m >$ is a presentation of the
fundamental group of $V$ in $v_0$ and $\Sigma$ is the analytic submanifold passing through the point
corresponding to $H$ and complementary to the orbit at this point.
Finally we get :

\begin{equation*}
    K^{h}_\F \simeq U_V\times K^{eq}_{\mathbb{CP}^1}.
\end{equation*}

Geometrically, this means that a holomorphic deformation corresponds to the data of two things. Firstly, a deformation of the complex structure of the base $V$ that induces a holomorphic deformation on all the leaves, while preserving the transverse structure of the foliation, given here by the fibration. Secondly, a deformation of the transverse structure, i.e. the fibration, by modifying the holonomy representation, we transform the way the leaves are transverse to the fibers.

\subsection*{iii) $\F$ is given as a fibration by hypersurfaces whose canonical bundle is trivial}

We denote $p:X\to C$ the corresponding fibration. For all $c\in C$, we have $p^{-1}(c)=H_c$ such that $c_1(H_c)= 0$.

The Leray spectral sequence gives exact sequences of low-degrees. For $\Theta^{tr}_\F$ we have :

\begin{equation*}
    0\xrightarrow[]{} H^1(C,p_*\Theta^{tr}_\F)\xrightarrow[]{} H^1(X,\Theta^{tr}_\F)\to H^0(C,\mathcal{R}^1p_*\Theta^{tr}_\F)\xrightarrow{} H^2(C,p_*\Theta^{tr}_\F) 
\end{equation*}

Moreover, we have $\Theta^{tr}_\F=p^{-1}TC$ and $\mathcal{R}^1 p_{*}\Theta^{tr}_\F=(\mathcal{R}^1 p_* \underline{\C})\otimes TC$ then :

\begin{equation*}
    0\xrightarrow[]{} H^1(C,TC)\xrightarrow[]{} H^1(X,\Theta^{tr}_\F)\to H^0(C,(\mathcal{R}^1 p_* \underline{\C})\otimes TC)\xrightarrow{} H^2(C,TC)=0.
\end{equation*}

When $H^0(C,(\mathcal{R}^1p_* \mathcal{O}_X)\otimes TC)=0$ (for example when $H_c$ is simply connected), we get $H^1(C,TC)\simeq H^1(X,\Theta^{tr}_\F)$ and then $K_C\simeq K^{tr}_\F$, where $K_C$ represents the Kuranishi space associated to the curve $C$ which is an open neighborhood $U_C$ of 0 in $H^1(C,TC)$ (since $H^2(C,TC)=0)$.\\
\\

To compute the $f$-deformations, we write the second exact sequence of low-degrees :

\begin{equation*}
    0\xrightarrow[]{} H^1(C,p_*T\F)\xrightarrow[]{} H^1(X,T\F)\to H^0(C,\mathcal{R}^1 p_*T\F)\xrightarrow{} H^2(C,p_* T\F)=0.
\end{equation*}

We recognize $ H^1(C,p_*T\F)$ the Čech $1$-cocyles of infinitesimal vertical automorphisms, i.e, the infinitesimal deformations corresponding to the $f$-deformations preserving the complex structure of leaves. We recognize also $H^0(C,\mathcal{R}^1 p_*T\F)$ as the holomorphic families of infinitesimal deformations of the $H_c$.

\pagebreak
\bibliographystyle{plain}
\bibliography{sample}

\end{document}